\documentclass[12pt]{article}
\usepackage{amssymb}
\usepackage{amsmath}
\usepackage{amscd}
\newtheorem{theorem}{Theorem}[section]
\newtheorem{proposition}[theorem]{Proposition}
\newtheorem{lemma}[theorem]{Lemma}
\newtheorem{corollary}[theorem]{Corollary}
\newtheorem{condition}[theorem]{Condition}
\newtheorem{definition}[theorem]{Definition}
\newtheorem{remark}[theorem]{Remark}

\newcommand{\bthm}{\begin{theorem}}
\newcommand{\ethm}{\end{theorem}}
\newcommand{\bpr}{\begin{proposition}}
\newcommand{\epr}{\end{proposition}}
\newcommand{\blem}{\begin{lemma}}
\newcommand{\elem}{\end{lemma}}
\newcommand{\bco}{\begin{corollary}}
\newcommand{\eco}{\end{corollary}}
\newcommand{\bcon}{\begin{condition}\rm}
\newcommand{\econ}{\end{condition}}
\newcommand{\bde}{\begin{definition}\rm}
\newcommand{\ede}{\end{definition}}
\newcommand{\bre}{\begin{remark}\rm}
\newcommand{\ere}{\end{remark}}
\newcommand{\bprf}{\noindent{\it Proof.\ }}
\newcommand{\eprf}{\hspace*{\fill} \rule{1.6mm}{3.2mm} \vspace{1.6mm}}
\newcommand{\benu}{\begin{enumerate}\renewcommand{\labelenumi}{{\rm (\roman{enumi})}}\renewcommand{\itemsep}{0pt}}
\newcommand{\eenu}{\end{enumerate}}

\newcommand{\N}{\mathbb{N}}
\newcommand{\Z}{\mathbb{Z}}
\newcommand{\R}{\mathbb{R}}
\newcommand{\C}{\mathbb{C}}
\newcommand{\T}{\mathbb{T}}
\newcommand{\e}{\varepsilon}
\newcommand{\K}{\mathbb{K}}
\newcommand{\M}{\mathbb{M}}

\DeclareMathOperator{\Prim}{Prim}
\DeclareMathOperator{\spa}{span}
\DeclareMathOperator{\cspa}{\overline{span}}
\newcommand{\On}{{\mathcal O}_n}
\newcommand{\Oi}{{\mathcal O}_\infty}
\newcommand{\cpi}{{\mathcal O}_\infty{\rtimes_{\alpha^\omega}}G}
\newcommand{\ip}[2]{\langle\,{#1}\,|\,{#2}\,\rangle}
\newcommand{\W}{{\mathcal W}}
\newcommand{\F}{{\mathcal F}}
\newcommand{\G}{{\mathcal G}}
\newcommand{\sg}{{\rm{sg}}(\omega)}
\newcommand{\csg}{\overline{\rm{sg}}(\omega)}
\newcommand{\vcsg}{\overline{\rm{sg}}_1(\omega)}

\begin{document}
\title{On crossed products of the Cuntz algebra ${\mathcal O}_\infty$ by quasi-free actions of abelian groups}
\author{Takeshi KATSURA\\
Department of Mathematical Sciences\\
University of Tokyo, Komaba, Tokyo, 153-8914, JAPAN\\
e-mail: {\tt katsu@ms.u-tokyo.ac.jp}}
\date{}

\maketitle

\begin{abstract}
{\footnotesize We investigate the structures of crossed products of the 
Cuntz algebra ${\mathcal O}_\infty$ by quasi-free actions of abelian groups.
We completely determine their ideal structures and compute the strong Connes
spectra and K-groups.}
\end{abstract}

\section{Introduction}

The crossed products of $C^*$-algebras give us plenty of interesting examples, 
and the structures of them have been examined by several authors.
In \cite{Ki}, A. Kishimoto gave a necessary and sufficient condition 
that the crossed products by abelian groups become simple 
in terms of the strong Connes spectrum.
For the case of the crossed products of Cuntz algebras 
by so-called quasi-free actions of abelian groups, 
he gave a condition for simplicity, which is easy to check.
In \cite{KK1} and \cite{KK2}, A. Kishimoto and A. Kumjian dealt with, 
among others, the crossed products of Cuntz algebras by quasi-free actions of 
the real group $\R$.
In our previous papers \cite{Ka1}, \cite{Ka2},
we examined the structures of crossed products of Cuntz algebras $\On$ 
by quasi-free actions of 
arbitrary locally compact, second countable, abelian groups.
The class of our algebras has many examples of simple stably projectionless 
$C^*$-algebras as well as AF-algebras and purely infinite $C^*$-algebras.
In \cite{Ka1}, we completely determined the ideal structures of our algebras, 
and gave another proof of A. Kishimoto's result on the simplicity of them. 
We also gave a necessary and sufficient condition that our algebras become 
primitive, and computed the Connes spectra and K-groups of our algebras.
In \cite{Ka2}, we proved that our algebras become AF-embeddable 
when actions satisfy certain conditions.
To the best of the author's knowledge, this is the first case 
to have succeeded in embedding crossed products of purely infinite 
$C^*$-algebras into AF-algebras except trivial cases.
We also gave a necessary and sufficient condition that our algebras become 
simple and purely infinite, 
and consequently our algebras are either purely infinite or AF-embeddable 
when they are simple.

In this paper, we deal with crossed products of the Cuntz algebra $\Oi$ 
by quasi-free actions of arbitrary locally compact, second countable, 
abelian groups.
From section 3 to section 6, we completely determine the ideal structures
of such algebras by using the technique developed in \cite{Ka1}.
We omit detailed computations if similar computations have been already done 
in \cite{Ka1}. 
Readers are referred to \cite{Ka1}.
In the last section, we gather some results on crossed products 
of the Cuntz algebra $\Oi$.
Among others, we give another proof of the determination of the simplicity 
of the crossed products done by A. Kishimoto, 
and we succeed in computing the strong Connes spectra of quasi-free actions 
on the Cuntz algebra $\Oi$.

The crossed products examined in this paper or in \cite{Ka1}, \cite{Ka2},
can be considered as continuous counterparts of 
Cuntz-Krieger algebras or graph algebras (cf. \cite{D}).
From this point of view, the crossed products of $\On$ can be considered 
as graph algebras of locally finite graphs, 
and the ones of $\Oi$ can be considered as graph algebras of graphs 
whose vertices emit and receive infinitely many edges.
Recently the ideal structures of graph algebras, which is not necessarily 
locally finite, were deeply examined in \cite{BHRS} and \cite{HS}.
Compared with row finite case,
it is rather difficult to describe ideal structures of graph algebras 
which have vertices emitting infinitely many edges.
This seems to be related to the difficulty of examination of the 
ideal structures of the crossed products of $\Oi$ 
compared with the ones of $\On$ done in \cite{Ka1}.

{\bf Acknowledgment.} 
The author would like to thank 
to his advisor Yasuyuki Kawahigashi for his support and encouragement, 
to Masaki Izumi for various comments and many suggestions.
He is also grateful to Iain Raeburn and Wojciech Szyma\'nski
for stimulating discussions. 
This work was partially supported by Research Fellowship 
for Young Scientists of the Japan Society for the Promotion of Science.

\section{Preliminaries}\label{PRE}

The Cuntz algebra $\Oi$ is the universal $C^*$-algebra generated 
by infinitely many isometries 
$S_1,S_2,\ldots$ satisfying $S_i^*S_j=\delta_{i,j}$.
For $n\in\Z_+:=\{1,2,\ldots\}$ and $k\in\N:=\{0,1,\ldots\}$, 
we define the set $\W_n^{(k)}$ of words in $\{1,2,\ldots,n\}$ 
with length $k$ by $\W_n^{(0)}=\{\emptyset\}$ and 
$$\W_n^{(k)}
=\big\{ (i_1,i_2,\ldots,i_k)\ \big|\ i_j\in\{1,2,\ldots,n\}\big\}$$
for $k\geq 1$.
Set $\W_n=\bigcup_{k=0}^\infty \W_n^{(k)}$ and 
$\W_\infty=\bigcup_{n=1}^\infty \W_n$.
For $\mu=(i_1,i_2,\ldots,i_k)\in\W_\infty$, 
we denote its length $k$ by $|\mu|$, 
and set $S_\mu=S_{i_1}S_{i_2}\cdots S_{i_k}\in\Oi$. 
Let $G$ be a locally compact abelian group which satisfies 
the second axiom of countability and $\Gamma$ be the dual group of $G$.
We use $+$ for multiplicative operations of abelian groups except for $\T$, 
which is the group of the unit circle in the complex plane $\C$. 
The pairing of $t\in G$ and $\gamma\in\Gamma$ is denoted by 
$\ip{t}{\gamma}\in\T$. 

For $\omega=(\omega_1,\omega_2,\ldots)\in\Gamma^\infty$,
we define an action $\alpha^\omega$ of abelian group $G$ on $\Oi$ by 
$\alpha^\omega_t(S_i)=\ip{t}{\omega_i}S_i$ for $i\in\Z_+$ and $t\in G$.
The action $\alpha^\omega:G\curvearrowright\Oi$ becomes quasi-free 
(for a definition of quasi-free actions on Cuntz algebras, see \cite{E}).
However, there exist quasi-free actions of abelian group $G$ on $\Oi$, 
which are not conjugate to $\alpha^\omega$ for any $\omega\in\Gamma^\infty$
though we do not deal with such actions.
The crossed product $\cpi$ has a $C^*$-subalgebra 
$\C 1{\rtimes_{\alpha^\omega}}G$ which is isomorphic to $C_0(\Gamma)$.
We consider $C_0(\Gamma)$ as a $C^*$-subalgebra of $\cpi$.
The Cuntz algebra $\Oi$ is naturally embedded into 
the multiplier algebra $M(\cpi)$ of $\cpi$. 
For each $\mu=(i_1,i_2,\ldots,i_k)\in\W_\infty$, 
we define an element $\omega_\mu$ of $\Gamma$ 
by $\omega_\mu=\sum_{j=1}^{k}\omega_{i_j}$.
For $\gamma_0\in\Gamma$, we define a (reverse) shift automorphism 
$\sigma_{\gamma_0}:C_0(\Gamma)\to C_0(\Gamma)$ by 
$(\sigma_{\gamma_0} f)(\gamma)=f(\gamma+\gamma_0)$ for $f\in C_0(\Gamma)$. 
Once noting that $\alpha^\omega_t(S_\mu)=\ip{t}{\omega_\mu}S_\mu$ 
for $\mu\in\W_\infty$, one can easily verify that 
$fS_\mu =S_\mu\sigma_{\omega_\mu}f$ for any $f\in C_0(\Gamma)\subset \cpi$.
For a subset $X$ of a $C^*$-algebra, 
we denote by $\spa X$ the linear span of $X$, 
and by $\cspa X$ its closure.
We have 
$\cpi=\cspa\{ S_\mu fS_\nu^*\mid \mu,\nu\in\W_\infty,\ f\in C_0(\Gamma)\}$.

We denote by $\M_k$ the $C^*$-algebra of $k \times k$ matrices 
for $k=1,2,\ldots$, 
and by $\K$ the $C^*$-algebra of compact operators 
of the infinite dimensional separable Hilbert space.

\section{Gauge invariant ideals}

In this section, we determine all the ideals 
which are globally invariant under the gauge action.
Here an ideal means a closed two-sided ideal, 
and the gauge action $\beta:\T\curvearrowright\cpi$ is defined by 
$\beta_t(S_\mu fS_\nu^*)=t^{|\mu|-|\nu|}S_\mu fS_\nu^*$ 
for $\mu,\nu\in\W_\infty,\ f\in C_0(\Gamma)$ and $t\in\T$.

For a positive integer $n$, we define a projection $p_n$ by 
$p_n=1-\sum_{i=1}^nS_iS_i^*$.
We set $p_0=1$.
Since $p_n$ commutes with $C_0(\Gamma)$, 
$p_n C_0(\Gamma)$ is a $C^*$-subalgebra of $\cpi$, 
which is isomorphic to $C_0(\Gamma)$.

\bde\label{omega}
Let $I$ be an ideal of the crossed product $\cpi$.
For each $n\in\N$, 
we define the closed subset $X_I^{(n)}$ of $\Gamma$ by
$$X_I^{(n)}=\{\gamma\in\Gamma\mid 
f(\gamma)=0\mbox{ for all }f\in C_0(\Gamma)\mbox{ with }p_nf\in I\}.$$
Set $X_I=X_I^{(0)}$, $X_I^{(\infty)}=\bigcap_{n=1}^\infty X_I^{(n)}$, 
and denote by $\widetilde{X}_I$ 
the pair $(X_I,X_I^{(\infty)})$ of subsets of $\Gamma$.
\ede

In other words, $X_I^{(n)}$ is determined by 
$p_n C_0(\Gamma\setminus X_I^{(n)})=I\cap p_n C_0(\Gamma)$.
One can easily see that 
$X_{I_1\cap I_2}^{(n)}=X_{I_1}^{(n)}\cup X_{I_2}^{(n)}$ for any $n\in\N$, 
hence $X_{I_1\cap I_2}=X_{I_1}\cup X_{I_2},\ 
X_{I_1\cap I_2}^{(\infty)}=X_{I_1}^{(\infty)}\cup X_{I_2}^{(\infty)}$ 
and that $I_1\subset I_2$ implies 
$X_{I_1}^{(n)}\supset X_{I_2}^{(n)}$ for any $n\in\N$, hence implies 
$X_{I_1}\supset X_{I_2},\ X_{I_1}^{(\infty)}\supset X_{I_2}^{(\infty)}$.
For $n\in\N$, the set $X_I^{(n)}$ can be described 
only in terms of $X_I$ and $X_I^{(\infty)}$.

\blem\label{X_I^{(n)}}
For an ideal $I$ of $\cpi$, we have
\begin{align*}
X_I^{(n)}&=X_I^{(\infty)}\cup\bigcup_{i=n+1}^\infty (X_I+\omega_i),
\end{align*}
for any $n\in\N$.
\elem

\bprf
Let $\gamma$ be an element of $X_I$ 
and $i$ be a positive integer grater than $n$.
Take $f\in C_0(\Gamma)$ with $p_nf\in I$.
Since 
$$S_{i}^*p_nfS_{i}=S_{i}^*fS_{i}=S_{i}^*S_{i}\sigma_{\omega_{i}}f
=\sigma_{\omega_{i}}f,$$ 
we have $\sigma_{\omega_{i}}f\in I\cap C_0(\Gamma)$.
Since $\gamma\in X_I$, we have $\sigma_{\omega_{i}}f(\gamma)=0$.
Hence $f(\gamma+\omega_{i})=0$ for any $f\in C_0(\Gamma)$ with $p_nf\in I$.
It implies $\gamma+\omega_{i}\in X_I^{(n)}$.
Thus $X_I^{(n)}\supset X_I+\omega_i$ for any $i>n$.
For $n\leq m$, we have $X_I^{(n)}\supset X_I^{(m)}$ because $p_np_m=p_m$. 
Therefore $X_I^{(n)}\supset X_I^{(\infty)}$.
Thus $X_I^{(n)}\supset 
X_I^{(\infty)}\cup\bigcup_{i=n+1}^\infty (X_I+\omega_i)$.

Conversely, take 
$\gamma\notin X_I^{(\infty)}\cup\bigcup_{i=n+1}^\infty (X_I+\omega_i)$.
Since $\gamma\notin X_I^{(\infty)}$, we can find a positive integer $m$ so that
$\gamma\notin X_I^{(m)}$.
When $m\leq n$, we see that $\gamma\notin X_I^{(n)}$.
We will show $\gamma\notin X_I^{(n)}$ in the case $m>n$.
Since $\gamma\notin X_I^{(m)}$, 
there exists $f\in C_0(\Gamma)$ such that $p_mf\in I$ and $f(\gamma)\neq 0$.
For each $i=n+1,n+2,\ldots,m$, there exists $f_i\in C_0(\Gamma)\cap I$ 
such that $f_i(\gamma-\omega_i)\neq 0$ because $\gamma\notin X_I+\omega_i$.
Set $g=f\prod_{i=n+1}^m\sigma_{-\omega_i}f_i$.
We have $g(\gamma)\neq 0$ and
$$p_n g=p_m g+\sum_{i=n+1}^mS_iS_i^*g=p_m g+\sum_{i=n+1}^mS_i(\sigma_{\omega_i}g)S_i^*\in I.$$
Therefore $\gamma\notin X_I^{(n)}$.
Thus we have 
$X_I^{(n)}=X_I^{(\infty)}\cup\bigcup_{i=n+1}^\infty (X_I+\omega_i)$.
\eprf

\bde
A subset $X$ of $\Gamma$ is called {\em $\omega$-invariant} 
if $X$ is a closed set with $X+\omega_i\subset X$ for any $i\in\Z_+$.
For an $\omega$-invariant set $X$, we define a closed set $H_X$ by
$$H_X=\overline{X\setminus\bigcup_{i=1}^\infty(X+\omega_i)}\ \cup\ 
\bigcap_{n=1}^\infty\overline{\bigcup_{i=n}^\infty(X+\omega_i)}.$$
\ede

Note that $H_X$ is a closed subset of $X$.

\bde
A pair $\widetilde{X}=(X,X^\infty)$ of subsets of $\Gamma$ is called 
{\em $\omega$-invariant} if $X$ is an $\omega$-invariant set, and 
$X^\infty$ is a closed set satisfying $H_X\subset X^\infty\subset X$.
\ede

\bpr
For any ideal $I$ of the crossed product $\cpi$, 
the pair $\widetilde{X}_I$ is $\omega$-invariant.
\epr

\bprf
By Lemma \ref{X_I^{(n)}}, 
we have $X_I=X_I^{(\infty)}\cup \bigcup_{i=1}^\infty(X_I+\omega_i)$.
From this, we see that $X_I$ is $\omega$-invariant and that 
$X_I\setminus\bigcup_{i=1}^\infty(X_I+\omega_i)\subset 
X_I^{(\infty)}\subset X_I$.
By Lemma \ref{X_I^{(n)}}, 
we have $\overline{\bigcup_{i=n}^\infty(X+\omega_i)}\subset\overline{X_I^{(n)}}
=X_I^{(n)}$.
Hence $\bigcap_{n=1}^\infty\overline{\bigcup_{i=n}^\infty(X+\omega_i)}\subset
\bigcap_{n=1}^\infty X_I^{(n)}=X_I^{(\infty)}$.
Therefore we get $H_X\subset X_I^{(\infty)}\subset X_I$.
\eprf

We will show that for an $\omega$-invariant pair $\widetilde{X}$, 
there exists a gauge invariant ideal $I$ such that 
$\widetilde{X}_I=\widetilde{X}$ (Proposition \ref{exist}).

\blem\label{X^{(n)}}
Let $\widetilde{X}=(X,X^{(\infty)})$ be an $\omega$-invariant pair. 
For $n\in\N$, set $X^{(n)}=X^{(\infty)}\cup\bigcup_{i=n+1}^\infty(X+\omega_i)$.
Then we have the following.
\benu
\item $X^{(n)}$ is closed for all $n\in\N$.
\item $X=X^{(0)}$, $X^{(\infty)}=\bigcap_{n=1}^\infty X^{(n)}$.
\item For $0\leq n<m$, $X^{(n)}=X^{(m)}\cup\bigcup_{i=n+1}^m(X+\omega_i)$.
\item For a positive integer $n$, 
$$X=\bigcup_{\mu\in\W_n}(X^{(n)}+\omega_\mu)\cup
\bigcap_{k=1}^\infty \bigg( \bigcup_{\mu\in\W_n^{(k)}}(X+\omega_\mu)\bigg).$$
\eenu
\elem

\bprf
\benu
\item Take $\gamma\in\overline{X^{(n)}}$ for a positive integer $n$.
If $U\cap X^{(\infty)}\neq\emptyset$ for all neighborhood $U$ of $\gamma$,
then $\gamma\in X^{(\infty)}\subset X^{(n)}$ because $X^{(\infty)}$ is closed.
Otherwise, we can find a positive integer $i_U$ grater than $n$ with 
$U\cap (X+\omega_{i_U})\neq\emptyset$ for any neighborhood $U$ of $\gamma$.
If there exists $i$ such that $i_U=i$ eventually, 
then $\gamma\in X+\omega_i\subset X^{(n)}$ because $X+\omega_i$ is closed.
If there are no such $i$, then we can see that 
$\gamma\in\overline{\bigcup_{i=m}^\infty(X+\omega_i)}$ for any $m$ with $m>n$.
Hence $\gamma\in H_{X}\subset X^{(\infty)}\subset X^{(n)}$.
Thus we have proved that $\gamma\in X^{(n)}$, 
from which it follows that $X^{(n)}$ is closed.
\item Since 
$X\setminus\bigcup_{i=1}^\infty(X+\omega_i)\subset X^{(\infty)}\subset X$, 
we have $X=X^{(0)}$.
We see that
\begin{align*}
\bigcap_{n=1}^\infty X^{(n)}
&=\bigcap_{n=1}^\infty
\bigg(X^{(\infty)}\cup\bigcup_{i=n+1}^\infty(X+\omega_i)\bigg)
=X^{(\infty)}\cup\bigcap_{n=1}^\infty
\bigg(\bigcup_{i=n+1}^\infty(X+\omega_i)\bigg).
\end{align*}
Since $\bigcap_{n=1}^\infty\left(\bigcup_{i=n+1}^\infty(X+\omega_i)\right)
\subset H_{X}\subset X^{(\infty)}$, 
we have $\bigcap_{n=1}^\infty X^{(n)}=X^{(\infty)}$.
\item It is obvious by the definition.
\item For a positive integer $n$, we have 
$X=X^{(n)}\cup\bigcup_{i=1}^n(X+\omega_i)$ by (iii).
Recursively, we get 
$X=\bigcup_{m=0}^{k-1}\big(\bigcup_{\mu\in\W_n^{(m)}}(X^{(n)}+\omega_\mu)\big)
\cup\bigcup_{\mu\in\W_n^{(k)}}(X+\omega_\mu)$ for any positive integer $k$.
Hence $X=\bigcup_{\mu\in\W_n}(X^{(n)}+\omega_\mu)\cup
\bigcap_{k=1}^\infty \big( \bigcup_{\mu\in\W_n^{(k)}}(X+\omega_\mu)\big)$.
\eprf
\eenu 

\bde
For an $\omega$-invariant pair $\widetilde{X}=(X,X^{(\infty)})$,
we define $I_{\widetilde{X}}\subset\cpi$ by
$$I_{\widetilde{X}}=\cspa\{S_\mu p_n fS_\nu^*\mid 
\mu,\nu\in\W_\infty,\ f\in C_0(\Gamma\setminus X^{(n)}),\ n\in\N\},$$
where $X^{(n)}=X^{(\infty)}\cup\bigcup_{i=n+1}^\infty(X+\omega_i)$.
\ede

\bpr\label{I_X,Xinfty}
For an $\omega$-invariant pair $\widetilde{X}=(X,X^{(\infty)})$,
the set $I_{\widetilde{X}}$ becomes a gauge invariant ideal of $\cpi$.
\epr

\bprf
Clearly $I_{\widetilde{X}}$ is a $*$-invariant closed linear space,
and is invariant under the gauge action $\beta$ because 
$\beta_t(S_\mu p_n fS_\nu^*)=t^{|\mu|-|\nu|}S_\mu p_n fS_\nu^*$ 
for $t\in\T$.
To prove that $I_{\widetilde{X}}$ is an ideal,
it suffices to show that for any $\mu_1,\nu_1,\mu_2,\nu_2\in\W_\infty$ 
and any $f\in C_0(\Gamma\setminus X^{(n)}),\ g\in C_0(\Gamma)$, 
the product $xy$ of $x=S_{\mu_1}p_n fS_{\nu_1}^*\in I_{\widetilde{X}}$ 
and $y=S_{\mu_2}gS_{\nu_2}^*\in\cpi$ is in $I_{\widetilde{X}}$.
If $S_{\nu_1}^*S_{\mu_2}=0$ or $S_{\nu_1}^*S_{\mu_2}=S_\mu^*$ 
for some $\mu\in\W_\infty$, 
then it is easy to see that $xy\in I_{\widetilde{X}}$.
Otherwise $S_{\nu_1}^*S_{\mu_2}=S_\mu$ 
for some $\mu=(i_1,i_2,\ldots,i_k)\in\W_\infty$ with $\mu\neq\emptyset$.
When $i_1\leq n$, we have $p_n fS_\mu=p_n S_\mu\sigma_{\omega_\mu}f=0$.
Hence $xy=0\in I_{\widetilde{X}}$.
When $i_1> n$, we have 
$p_n fS_\mu=p_n S_\mu\sigma_{\omega_\mu}f=S_\mu\sigma_{\omega_\mu}f$.
Now, $f\in C_0(\Gamma\setminus X^{(n)})$ implies 
$\sigma_{\omega_\mu}f\in C_0(\Gamma\setminus X)$ because 
$X+\omega_\mu\subset X+\omega_{i_1}\subset X^{(n)}$.
Hence we have $xy\in I_{\widetilde{X}}$.
It completes the proof.
\eprf

\bpr\label{exist}
Let $\widetilde{X}=(X,X^{(\infty)})$ be an $\omega$-invariant pair,
and set $I=I_{\widetilde{X}}$. 
Then $\widetilde{X}_I=\widetilde{X}$.
\epr

\bprf
By the definition of $I$, we get $X_I^{(n)}\subset X^{(n)}$ for any $n\in\N$.
We will first prove that $X_I=X$.
To the contrary, assume that $X_I\subsetneqq X$.
Then there exists $f\in I\cap C_0(\Gamma)$ 
such that $f(\gamma_0)=1$ for some $\gamma_0\in X$.
Since $f\in I$, there exist $n_l\in\N$, 
$f_l\in C_0(\Gamma\setminus X^{(n_l)})$ and 
$\mu_l,\nu_l\in\W_\infty\ (l=1,2,\ldots,L)$ such that 
$$\bigg\| f-\sum_{l=1}^LS_{\mu_l}p_{n_l}f_lS_{\nu_l}^*\bigg\|<\frac12.$$
Take a positive integer $n$ so large that 
$n_l\leq n$ and $\mu_l,\nu_l\in\W_n$ for $l=1,2,\ldots,L$.
For any $\mu_0\in\W_n$, 
we have $p_nS_{\mu_0}^*fS_{\mu_0}p_n=p_n\sigma_{\omega_{\mu_0}}f$ 
and $\sigma_{\omega_{\mu_0}}f(\gamma_0-\omega_{\mu_0})=1$.
For $l$ with $\mu_l=\nu_l=\mu_0$, we have 
$p_nS_{\mu_0}^*(S_{\mu_l}p_{n_l}f_lS_{\nu_l}^*)S_{\mu_0} p_n=p_n f_l$.
For $l$ with $\mu_l\nu=\nu_l\nu=\mu_0$ for some 
$\nu=(i_1,i_2,\ldots,i_k)\in\W_n$ with $i_1>n_l$,
we have 
$p_nS_{\mu_0}^*(S_{\mu_l}p_{n_l}f_lS_{\nu_l}^*)S_{\mu_0} p_n
=p_n\sigma_{\omega_{\nu}}f_l$.
We have $\sigma_{\omega_{\nu}}f_l\in C_0(\Gamma\setminus X)$, 
because $X+\omega_{\nu}\subset X+\omega_{i_1}\subset X^{(n_l)}$.
For other $l$, we have 
$p_nS_{\mu_0}^*(S_{\mu_l}p_{n_l}f_lS_{\nu_l}^*)S_{\mu_0} p_n=0$.
Hence we get
$$\bigg\|\sigma_{\omega_{\mu_0}}f-\sum_{l=1}^L g_l\bigg\|=
\bigg\|p_n\bigg(\sigma_{\omega_{\mu_0}}f-\sum_{l=1}^L g_l\bigg)\bigg\|=\bigg\|p_nS_{\mu_0}^*\bigg(f-\sum_{l=1}^LS_{\mu_l}p_{n_l}f_lS_{\nu_l}^*\bigg)
S_{\mu_0}p_n\bigg\|<\frac12,$$
where $g_l\in C_0(\Gamma\setminus X^{(n_l)})$ when $\mu_l=\nu_l=\mu_0$, 
and $g_l\in C_0(\Gamma\setminus X)$ when $\mu_l\nu=\nu_l\nu=\mu_0$ for some 
$\nu=(i_1,i_2,\ldots,i_k)\in\W_n$ with $i_1>n_l$, and $g_l=0$ otherwise.
To derive a contradiction, it suffices to find $\mu_0\in\W_n$ such that
$g_l(\gamma_0-\omega_{\mu_0})=0$ for any $l$.
By Lemma \ref{X^{(n)}} (iv), we have either $\gamma_0\in \bigcap_{m=1}^\infty 
\big( \bigcup_{\mu\in\W_n^{(m)}}(X+\omega_\mu)\big)$ 
or $\gamma_0\in X^{(n)}+\omega_\mu$ for some $\mu\in\W_n$.

When $\gamma_0\in \bigcap_{m=1}^\infty 
\big( \bigcup_{\mu\in\W_n^{(m)}}(X+\omega_\mu)\big)$, 
take $\mu_0\in\W_n$ so that $|\mu_0|>|\mu_l|,|\nu_l|$ for $l=1,2,\ldots,L$ 
and $\gamma_0\in X+\omega_{\mu_0}$.
Then $\mu_l=\nu_l=\mu_0$ never occurs. 
Hence $g_l\in C_0(\Gamma\setminus X)$ for any $l$.
We get $g_l(\gamma_0-\omega_{\mu_0})=0$ because $\gamma_0-\omega_{\mu_0}\in X$.
When $\gamma_0\in X^{(n)}+\omega_\mu$ for some $\mu\in\W_n$, take $\mu_0=\mu$.
Since $\gamma_0-\omega_{\mu_0}\in X^{(n)}\subset X^{(n_l)}\subset X$,
we have $g_l(\gamma_0-\omega_{\mu_0})=0$ 
either if $g_l\in C_0(\Gamma\setminus X^{(n_l)})$ or 
if $g_l\in C_0(\Gamma\setminus X)$.
Hence $g_l(\gamma_0-\omega_{\mu_0})=0$ for any $l$.
Therefore we have $X_I=X$.

Next we will show that $X_I^{(n)}=X^{(n)}$ for a positive integer $n$.
To derive a contradiction, assume that $X_I^{(n)}\subsetneqq X^{(n)}$.
Then there exists $f\in C_0(\Gamma)$ 
such that $p_nf\in I$ and $f(\gamma_0)=1$ for some $\gamma_0\in X^{(n)}$.
Since $p_nf\in I$, there exist $n_l\in\N$, 
$f_l\in C_0(\Gamma\setminus X^{(n_l)})$ and 
$\mu_l,\nu_l\in\W_\infty\ (l=1,2,\ldots,L)$ such that 
$$\bigg\| p_nf-\sum_{l=1}^LS_{\mu_l}p_{n_l}f_lS_{\nu_l}^*\bigg\|<\frac12.$$
Take a positive integer $m$ so large that $\mu_l,\nu_l\in\W_{m}$, $n_l\leq m$ 
for $l=1,2,\ldots,L$ and $n\leq m$.
By Lemma \ref{X^{(n)}} (iii), we have 
$X^{(n)}=X^{(m)}\cup\bigcup_{i=n+1}^{m}(X+\omega_i)$.
When $\gamma_0\in X^{(m)}$, we have $f_l(\gamma_0)=0$ for any $l$.
On the other hand, we get $\|f-\sum_{\mu_l=\nu_l=\emptyset}f_l\|<1/2$
because
$$p_{m}\bigg(p_nf-\sum_{l=1}^LS_{\mu_l}p_{n_l}f_lS_{\nu_l}^*
\bigg)p_{m}=p_{m}f-\sum_{\mu_l=\nu_l=\emptyset}p_{m}f_l.$$
This is a contradiction.
When $\gamma_0\in X+\omega_i$ for some $i$ with $n<i\leq m$, 
we have $\sigma_{\omega_i}f=S_i^*(p_nf)S_i\in I$ and 
$\sigma_{\omega_i}f(\gamma_0-\omega_i)=1$.
This contradicts the fact that $X_I=X$.
Therefore $X_I^{(n)}=X^{(n)}$ for a positive integer $n$.
Hence $X_I^{(\infty)}=\bigcap_{n=1}^\infty X_I^{(n)}=
\bigcap_{n=1}^\infty X^{(n)}=X^{(\infty)}$.
We have shown that $\widetilde{X}_I=\widetilde{X}$.
\eprf

By Proposition \ref{exist}, the map $I\mapsto \widetilde{X}_I$ 
from the set of gauge invariant ideals $I$ of $\cpi$ 
to the set of $\omega$-invariant pairs is surjective.
Now, we turn to showing that this map is injective (Proposition \ref{unique}).
To do so, we investigate the quotient $(\cpi)/I$ of $\cpi$ by an ideal $I$ 
which is not $\cpi$. 
Since $I\cap C_0(\Gamma)=C_0(\Gamma\setminus X_I)$, a $C^*$-subalgebra 
$C_0(\Gamma)/(I\cap C_0(\Gamma))$ of $(\cpi)/I$ is isomorphic to $C_0(X_I)$.
We will consider $C_0(X_I)$ as a $C^*$-subalgebra of $(\cpi)/I$.
We will use the same symbols
$S_1,S_2,\ldots\in M((\cpi)/I)$ as the ones in $M(\cpi)$ for denoting 
the isometries of $\Oi$ which is naturally embedded into $M((\cpi)/I)$.
For an $\omega$-invariant set $X$, 
we can define a $*$-homomorphism $\sigma_{\omega_\mu}:C_0(X)\to C_0(X)$ 
for $\mu\in\W_\infty$.
This map $\sigma_{\omega_\mu}$ is always surjective, 
but it is injective only in the case that $X\subset X+\omega_\mu$, 
which is equivalent to $X=X+\omega_\mu$.
One can easily verify the following.

\blem\label{cp/I}
Let $I$ be an ideal that is not $\cpi$. 
For $\mu,\nu\in\W_\infty$ and $f\in C_0(X_I)\subset(\cpi)/I$, 
the following hold.
\benu
\item $S_\mu fS_\nu^*=0$ if and only if $f=0$.
\item For $n\in\N$, $p_nf=0$ if and only if $f\in C_0(X_I\setminus X_I^{(n)})$.
\item $fS_\mu=S_\mu \sigma_{\omega_\mu}f$.
\item $(\cpi)/I=\cspa\{S_\mu fS_\nu^*\mid\mu,\nu\in\W_\infty,\ f\in C_0(X_I)\}$.
\eenu
\elem

We define a $C^*$-subalgebra of $(\cpi)/I$, 
which corresponds to the AF-core for Cuntz algebras.

\bde\label{AFcore}
Let $I$ be an ideal that is not $\cpi$.
We define $C^*$-subalgebras of $(\cpi)/I$ by
\begin{align*}
\G_I^{(n,k)}&=\spa\{S_\mu fS_\nu^*\mid \mu,\nu\in\W_n^{(k)},\ 
f\in C_0(X_I)\},\\
\F_I^{(n,k)}&=\spa\{S_\mu p_nfS_\nu^*\mid \mu,\nu\in\W_n^{(k)},\ 
f\in C_0(X_I)\},\\
\F_I^{(n)}&=\spa\{S_\mu fS_\nu^*\mid \mu,\nu\in\W_n,\ 
0\leq |\mu|=|\nu|\leq n,\ f\in C_0(X_I)\},\\
\F_I&=\cspa\{S_\mu fS_\nu^*\mid \mu,\nu\in\W_\infty,\ 
|\mu|=|\nu|,\ f\in C_0(X_I)\},
\end{align*}
for $n\in\Z_+, 0\leq k\leq  n$.
\ede

\blem\label{F_n}
Let $I$ be an ideal that is not $\cpi$. 
For $n\in\Z_+, 0\leq k\leq  n$, we have the following.
\benu
\item $\G_I^{(n,k)}\cong \M_{n^k}\otimes C_0(X_I)$.
\item $\F_I^{(n,k)}\cong \M_{n^k}\otimes C_0(X_I^{(n)})$.
\item $\F_I^{(n)}\cong\bigoplus_{k=0}^{n-1}\F_I^{(n,k)}\oplus \G_I^{(n,n)}$.
\item $\bigcup_{n=1}^\infty\F_I^{(n)}$ is dense in $\F_I$.
\eenu
\elem

\bprf
\benu
\item Since the set $\W_n^{(k)}$ has $n^k$ elements, we may use 
$\{e_{\mu,\nu}\}_{\mu,\nu\in\W_n^{(k)}}$ for denoting the matrix units of 
$\M_{n^k}$.
One can easily see that 
$$\M_{n^k}\otimes C_0(X_I)\ni e_{\mu,\nu}\otimes f 
\mapsto S_\mu fS_\nu^*\in \G_I^{(n,k)}$$
gives us an isomorphism from $\M_{n^k}\otimes C_0(X_I)$ to $\G_I^{(n,k)}$.
\item We can define a surjective map from $\G_I^{(n,k)}$ to $\F_I^{(n,k)}$ by
$$\G_I^{(n,k)}\ni S_\mu fS_\nu^*\mapsto S_\mu p_n fS_\nu^*\in \F_I^{(n,k)}.$$
Its kernel is $\M_{n^k}\otimes C_0(X_I\setminus X_I^{(n)})$ 
under the isomorphism $\G_I^{(n,k)}\cong \M_{n^k}\otimes C_0(X_I)$ 
by Lemma \ref{cp/I} (ii).
Hence we have $\F_I^{(n,k)}\cong \M_{n^k}\otimes C_0(X_I^{(n)})$.
\item It can be done just by computation.
\item Obvious by the definitions of $\F_I^{(n)}$ and $\F_I$.
\eprf
\eenu

We will often identify $\G_I^{(n,n)}$ with $C_0(X_I,\M_{n^n})$. 
The following lemma essentially appeared in \cite{C}.

\blem\label{isom}
For $i=1,2$, let $E_i$ be a conditional expectation 
from a $C^*$-algebra $A_i$ onto a $C^*$-subalgebra $B_i$ of $A_i$.
Let $\varphi:A_1\to A_2$ be a $*$-homomorphism 
with $\varphi\circ E_1=E_2\circ\varphi$.
If the restriction of $\varphi$ on $B_1$ is injective and $E_1$ is faithful, 
then $\varphi$ is injective.
\elem

For an ideal $I$ which is invariant under the gauge action $\beta$, 
we can extend the gauge action on $\cpi$ to one on $(\cpi)/I$, 
which is also denoted by $\beta$.
The following lemma is standard.

\blem\label{cond.exp1}
Let $I$ be a gauge invariant ideal that is not $\cpi$. 
Then, 
$$E_{I}:(\cpi)/I\ni x\mapsto\int_{\T}\beta_t(x)dt\in (\cpi)/I$$
is a faithful conditional expectation onto $\F_{I}$, 
where $dt$ is the normalized Haar measure on $\T$.
\elem

\bpr\label{unique}
For any gauge invariant ideal $I$, we have $I_{\widetilde{X}_I}=I$.
\epr

\bprf
When $I=\cpi$, we have $X_I=X_I^{(\infty)}=\emptyset$. 
Thus $I_{\widetilde{X}_I}=\cpi$.
Let $I$ be a gauge invariant ideal that is not $\cpi$ 
and set $J=I_{\widetilde{X}_I}$. 
By the definition, $J\subset I$. 
Hence there exists a surjective $*$-homomorphism $\pi:(\cpi)/J\to(\cpi)/I$. 
By Proposition \ref{exist} and  Lemma \ref{F_n}, 
the restriction of $\pi$ on $\F_{J}^{(k)}$ 
is an isomorphism from $\F_{J}^{(k)}$ onto $\F_{I}^{(k)}$ 
and so the restriction of $\pi$ on $\F_{J}$ is an isomorphism 
from $\F_{J}$ onto $\F_{I}$.
By Lemma \ref{cond.exp1}, there are faithful conditional expectations 
$E_J:(\cpi)/J\to \F_{J}$ and $E_I:(\cpi)/I\to\F_{I}$
with $E_I\circ \pi=\pi\circ E_{J}$.
By Lemma \ref{isom}, $\pi$ is injective.
Therefore $I_{\widetilde{X}_I}=I$.
\eprf

\bthm\label{OneToOne}
The maps $I\mapsto \widetilde{X}_I$ and 
$\widetilde{X}\mapsto I_{\widetilde{X}}$ induce 
a one-to-one correspondence 
between the set of gauge invariant ideals of $\cpi$ 
and the set of $\omega$-invariant pairs of subsets of $\Gamma$.
\ethm

\bprf
Combine Proposition \ref{exist} and Proposition \ref{unique}.
\eprf

\section{Primeness for $\omega$-invariant pairs}

In this section, we give a necessary condition for an ideal to be primitive
in terms of $\omega$-invariant pairs.
We will use it after in order to determine all primitive ideals.

An ideal of a $C^*$-algebra is called primitive 
if it is a kernel of some irreducible representation.
A $C^*$-algebra is called primitive if $0$ is a primitive ideal.
When a $C^*$-algebra $A$ is separable, an ideal $I$ of $A$ is primitive 
if and only if $I$ is prime, i.e. 
for two ideals $I_1,I_2$ of $A$,
$I_1\cap I_2\subset I$ implies either $I_1\subset I$ or $I_2\subset I$.
We define primeness for $\omega$-invariant pairs.
For two $\omega$-invariant pair $\widetilde{X}_1=(X_1,X_1^{(\infty)}),\ 
\widetilde{X}_2=(X_2,X_2^{(\infty)})$, 
we write $\widetilde{X}_1\subset \widetilde{X}_2$ 
if $X_1\subset X_2, X_1^{(\infty)}\subset X_2^{(\infty)}$
and denote by $\widetilde{X}_1\cup\widetilde{X}_2$ 
the $\omega$-invariant pair $(X_1\cup X_2,X_1^{(\infty)}\cup X_2^{(\infty)})$.

\bde
An $\omega$-invariant pair $\widetilde{X}$ is called {\em prime} 
if $\widetilde{X}_1\cup \widetilde{X}_2\supset\widetilde{X}$ 
implies either $\widetilde{X}_1\supset\widetilde{X}$ 
or $\widetilde{X}_2\supset\widetilde{X}$
for two $\omega$-invariant pairs 
$\widetilde{X}_1,\ \widetilde{X}_2$.
\ede

\bpr\label{prime}
If an ideal $I$ of $\cpi$ is primitive, 
then $\widetilde{X}_I$ is a prime $\omega$-invariant pair.
\epr

\bprf
Let $I$ be a primitive ideal of $\cpi$.
Take two $\omega$-invariant pairs 
$\widetilde{X}_1$, $\widetilde{X}_2$ with
$\widetilde{X}_1\cup \widetilde{X}_2\supset \widetilde{X}_I$.
Set $I_1=I_{\widetilde{X}_1}$ and $I_2=I_{\widetilde{X}_2}$.
Then 
$$I_1\cap I_2=I_{\widetilde{X}_1\cup \widetilde{X}_2}\subset 
I_{\widetilde{X}_I}\subset I.$$
Since $I$ is prime, we have either $I_1\subset I$ or $I_2\subset I$.
Hence we get either $\widetilde{X}_1\supset \widetilde{X}_I$ 
or $\widetilde{X}_2\supset\widetilde{X}_I$.
Thus $\widetilde{X}_I$ is prime.
\eprf

In general, the converse of Proposition \ref{prime} is not true 
(see Corollary \ref{primitive1} and Proposition \ref{primitive2}).
The ideal $I$ is prime if and only if 
the equality $I_1\cap I_2=I$ implies either $I_1=I$ or $I_2=I$ 
for two ideals $I_1,I_2$
(see the proof of (iii)$\Rightarrow$(iv) of Proposition \ref{primepair}).
The following is the counterpart of this fact 
for prime $\omega$-invariant pairs.

\bpr\label{primepair}
For an $\omega$-invariant pair $\widetilde{X}$, 
the following are equivalent.
\benu
\item $\widetilde{X}$ is prime.
\item 
For two $\omega$-invariant pairs $\widetilde{X}_1$, 
$\widetilde{X}_2$,
the equality $\widetilde{X}_1\cup \widetilde{X}_2=\widetilde{X}$
implies either $\widetilde{X}_1=\widetilde{X}$ 
or $\widetilde{X}_2=\widetilde{X}$.
\item For two gauge invariant ideals $I_1,I_2$ of $\cpi$,
the equality $I_1\cap I_2=I_{\widetilde{X}}$ implies either
$I_1=I_{\widetilde{X}}$ or $I_2=I_{\widetilde{X}}$.
\item For two gauge invariant ideals $I_1,I_2$ of $\cpi$, 
the inclusion $I_1\cap I_2\subset I_{\widetilde{X}}$ implies 
either $I_1\subset I_{\widetilde{X}}$ or $I_2\subset I_{\widetilde{X}}$.
\eenu
\epr

\bprf
(i)$\Rightarrow$(ii): 
Take two $\omega$-invariant pairs $\widetilde{X}_1$, 
$\widetilde{X}_2$ with $\widetilde{X}_1\cup \widetilde{X}_2=\widetilde{X}$.
By (i), we have either $\widetilde{X}_1\supset \widetilde{X}$ 
or $\widetilde{X}_2\supset \widetilde{X}$.
Hence we get either $\widetilde{X}_1=\widetilde{X}$ 
or $\widetilde{X}_2=\widetilde{X}$.

(ii)$\Rightarrow$(iii): 
Take two gauge invariant ideals $I_1,I_2$
with $I_1\cap I_2=I_{\widetilde{X}}$.
We have $\widetilde{X}_{I_1}\cup \widetilde{X}_{I_2}=\widetilde{X}$.
By (ii), we have either $\widetilde{X}_{I_1}=\widetilde{X}$ or
$\widetilde{X}_{I_2}=\widetilde{X}$.
By Proposition \ref{unique},
we have either $I_1=I_{\widetilde{X}}$ or $I_2=I_{\widetilde{X}}$.

(iii)$\Rightarrow$(iv): 
Take two gauge invariant ideals $I_1,I_2$
with $I_1\cap I_2\subset I_{\widetilde{X}}$.
Then we have 
$$(I_1+I_{\widetilde{X}})\cap (I_2+I_{\widetilde{X}})=
(I_1\cap I_2) +I_{\widetilde{X}}=I_{\widetilde{X}}.$$
By (iii), either $I_1+I_{\widetilde{X}}=I_{\widetilde{X}}$ or 
$I_2+I_{\widetilde{X}}=I_{\widetilde{X}}$ holds.
Hence we get 
either $I_1\subset I_{\widetilde{X}}$ or $I_2\subset I_{\widetilde{X}}$.

(iv)$\Rightarrow$(i): 
Similarly as the proof of Proposition \ref{prime}.
\eprf

We will use the implication (ii)$\Rightarrow$(i) 
to determine which $\omega$-invariant pair is prime.
We also need a notion of primeness for $\omega$-invariant sets.

\bde
An $\omega$-invariant set $X$ is called {\em prime} 
if $X_1\cup X_2\supset X$ implies either $X_1\supset X$ or $X_2\supset X$,
for any $\omega$-invariant sets $X_1,X_2$.
\ede

We set $\sg=\{\omega_\mu\mid \mu\in\W_\infty\}$ 
which is the semigroup generated by $\omega_1,\omega_2,\ldots$ and 
denote by $\csg$ its closure.
Note that a closed subset $X$ of $\Gamma$ is $\omega$-invariant 
if and only if $X+\csg=X$.
For any $\gamma\in\Gamma$, it is easy to see that 
the set $\gamma+\csg$ is a prime $\omega$-invariant set.
The following is a necessary and sufficient condition 
for an $\omega$-invariant set to be prime, 
which can be considered as an analogue of maximal tails in \cite{BHRS}.

\bpr\label{Xprime}
An $\omega$-invariant set $X$ of $\Gamma$ is prime if and only if for any 
$\gamma_1,\gamma_2\in X$ and 
any neighborhoods $U_1$, $U_2$ of $\gamma_1,\gamma_2$ respectively,
there exist $\gamma\in X$ and $\mu_1,\mu_2\in\W_\infty$ with 
$\gamma+\omega_{\mu_1}\in U_1$ and $\gamma+\omega_{\mu_2}\in U_2$.
\epr

\bprf
Suppose $X$ is a prime $\omega$-invariant set.
Take $\gamma_1,\gamma_2\in X$ and 
neighborhoods $U_1$, $U_2$ of $\gamma_1,\gamma_2$ respectively.
Set $X_j=\Gamma\setminus\bigcup_{\mu\in\W_\infty}(U_j-\omega_\mu)$ for $j=1,2$.
Then $X_1$ and $X_2$ are $\omega$-invariant sets satisfying 
$X_1\not\supset X$ and $X_2\not\supset X$. 
Since $X$ is prime, we have $X_1\cup X_2\not\supset X$.
Hence there exists $\gamma\in X$ with $\gamma\notin X_1\cup X_2$.
By the definition of $X_1$ and $X_2$, there exist $\mu_1,\mu_2$ such that
$\gamma+\omega_{\mu_1}\in U_1$ and $\gamma+\omega_{\mu_2}\in U_2$.

Conversely assume that for any $\gamma_1,\gamma_2\in X$ and 
any neighborhoods $U_1$, $U_2$ of $\gamma_1,\gamma_2$ respectively,
there exist $\gamma\in X$ and $\mu_1,\mu_2\in\W_\infty$ with 
$\gamma+\omega_{\mu_1}\in U_1$ and $\gamma+\omega_{\mu_2}\in U_2$.
Take $\omega$-invariant sets $X_1$ and $X_2$ satisfying 
$X_1\not\supset X$ and $X_2\not\supset X$.
There exist $\gamma_1,\gamma_2\in X$ 
with $\gamma_1\notin X_1$ and $\gamma_2\notin X_2$.
Hence there exist $\gamma\in X$ and $\mu_1,\mu_2\in\W_\infty$ with 
$\gamma+\omega_{\mu_1}\notin X_1$ and $\gamma+\omega_{\mu_2}\notin X_2$.
Since $X_1$ and $X_2$ are $\omega$-invariant, we have $\gamma\notin X_1$
and $\gamma\notin X_2$.
Therefore, $X_1\cup X_2\not\supset X$.
Thus, $X$ is prime.
\eprf

\blem\label{primepair0}
If an $\omega$-invariant pair $\widetilde{X}=(X,X^{(\infty)})$ is prime, 
then $X^{(\infty)}=H_X$ 
or $X^{(\infty)}=H_X\cup\{\gamma\}$ for some $\gamma\notin H_X$.
\elem

\bprf
Let $\widetilde{X}=(X,X^{(\infty)})$ be a prime $\omega$-invariant pair.
To derive a contradiction, assume $X^{(\infty)}\setminus H_X$ has 
two points $\gamma_1,\gamma_2$.
Take open sets $U_1\ni\gamma_1$, $U_2\ni\gamma_2$ with $U_1\cap U_2=\emptyset$,
$U_1\cap H_X=\emptyset$ and $U_2\cap H_X=\emptyset$.
Then $\widetilde{X}_i=(X,X^{(\infty)}\setminus U_i)\ (i=1,2)$ 
are $\omega$-invariant pairs satisfying 
$\widetilde{X}=\widetilde{X}_1\cup \widetilde{X}_2$.
However, we have $\widetilde{X}\not\subset \widetilde{X}_1$ and 
$\widetilde{X}\not\subset \widetilde{X}_2$.
This contradicts the primeness of $\widetilde{X}$.
\eprf

\blem\label{primepair1}
An $\omega$-invariant pair $(X,H_X)$ is prime if and only if
$X$ is a prime $\omega$-invariant set.
\elem

\bprf
Suppose that $(X,H_X)$ is a prime $\omega$-invariant pair.
Take $\omega$-invariant sets $X_1,X_2$ with $X\subset X_1\cup X_2$.
We have $(X,H_X)\subset (X_1,X_1)\cup (X_2,X_2)$.
Since $(X,H_X)$ is prime, 
either $(X,H_X)\subset (X_1,X_1)$ or $(X,H_X)\subset (X_2,X_2)$ holds.
Therefore $X$ is a prime $\omega$-invariant set.
Conversely assume that $X$ is a prime $\omega$-invariant set.
Take two $\omega$-invariant pairs $(X_1,X_1^{(\infty)})$, 
$(X_2,X_2^{(\infty)})$ with 
$(X_1,X_1^{(\infty)})\cup (X_2,X_2^{(\infty)})=(X,H_X)$.
Since $X$ is prime, either $X\subset X_1$ or $X\subset X_2$.
We may assume $X\subset X_1$.
Then $X=X_1$.
Hence $H_X=H_{X_1}\subset X_1^{(\infty)}\subset  H_X$.
We get $(X_1,X_1^{(\infty)})=(X,H_X)$.
By Proposition \ref{primepair}, $(X,H_X)$ is a prime $\omega$-invariant pair.
\eprf

\blem\label{primepair2}
An $\omega$-invariant pair $(X,H_X\cup\{\gamma\})$ is prime 
for some $\gamma\notin H_X$ if and only if $X=\gamma+\csg$.
\elem

\bprf
Suppose that an $\omega$-invariant pair $(X,H_X\cup\{\gamma\})$ is prime.
Then $(X,H_X\cup\{\gamma\})\subset (X,H_X)\cup (\gamma+\csg,\gamma+\csg)$
implies $(X,H_X\cup\{\gamma\})\subset(\gamma+\csg,\gamma+\csg)$
because $H_X\cup\{\gamma\}\not\subset H_X$.
Hence $\gamma+\csg\subset X\subset\gamma+\csg$.
Thus, we get $X=\gamma+\csg$.
Conversely, assume $X=\gamma+\csg$.
Take two $\omega$-invariant pairs $(X_1,X_1^{(\infty)})$, 
$(X_2,X_2^{(\infty)})$ with 
$(X_1,X_1^{(\infty)})\cup (X_2,X_2^{(\infty)})=(X,H_X\cup\{\gamma\})$.
We may assume $\gamma\in X_1^{(\infty)}$.
Then we have $X=\gamma+\csg\subset X_1^{(\infty)}+\csg\subset X_1\subset X$.
Hence $X_1=X$.
We have $H_X\cup\{\gamma\}=H_{X_1}\cup\{\gamma\}\subset 
X_1^{(\infty)}\subset H_X\cup\{\gamma\}$.
Therefore $(X_1,X_1^{(\infty)})=(X,H_X\cup\{\gamma\})$.
By Proposition \ref{primepair}, $(X,H_X\cup\{\gamma\})$ is a prime $\omega$-invariant pair.
\eprf

\bpr\label{pp}
An $\omega$-invariant pair $(X,X^{(\infty)})$ is prime if and only if
either $X$ is prime and $X^{(\infty)}=H_X$ or $X=\gamma+\csg$ and 
$X^{(\infty)}=H_X\cup\{\gamma\}$ for some $\gamma\notin H_X$.
\epr

\bprf
Combine Lemma \ref{primepair0}, Lemma \ref{primepair1} 
and Lemma \ref{primepair2}.
\eprf

\section{The ideal structure of $\cpi$ (part 1)}

In this section and the next section, 
we completely determine the ideal structure of $\cpi$ 
(Theorem \ref{idestr1}, Theorem \ref{idestr2}).
The ideal structure of $\cpi$ depends on whether $\omega\in\Gamma^\infty$ satisfies 
the following condition:

\bcon\label{cond}
For each $i\in\Z_+$, one of the following two conditions is satisfied:
\benu
\item For any positive integer $k$, $k\omega_i\neq 0$.
\item There exists a sequence $\mu_1,\mu_2,\ldots$ in $\W_\infty$ such that
$S_{\mu_k}^*S_i=0$ for any $k$ and $\lim_{k\to\infty}\omega_{\mu_k}=0$.
\eenu
\econ

This condition is an analogue of Condition (K) in the case of 
graph algebras \cite{BHRS}.
In this section, 
we deal with the case that $\omega$ satisfies Condition \ref{cond}.

\bpr\label{cond.exp2}
If $\omega$ satisfies Condition \ref{cond},
then for an ideal $I$ that is not $\cpi$,
there exists a unique conditional expectation $E_I$ from $(\cpi)/I$ onto 
$\F_I$ such that 
$E_I(S_\mu fS_\nu^*)=\delta_{|\mu|,|\nu|}S_\mu fS_\nu^*$ 
for $\mu,\nu\in\W_\infty,\ f\in C_0(X_I)$.
\epr

\bprf
Take $x=\sum_{l=1}^LS_{\mu_l}f_lS_{\nu_l}^*\in (\cpi)/I$ where
$\mu_l,\nu_l\in\W_\infty$ and $f_l\in C_0(X_I)$ for $l=1,2,\ldots,L$.
Set $x_0=\sum_{|\mu_l|=|\nu_l|}S_{\mu_l}f_lS_{\nu_l}^*$ 
and we will prove that $\|x_0\|\leq\|x\|$.
If we choose a positive integer $n$ so that $|\mu_l|,|\nu_l|\leq n$ and 
$\mu_l,\nu_l\in\W_n$ for $l=1,2,\ldots,L$,
then $x_0\in\F_I^{(n)}$.
By Lemma \ref{F_n}, 
there exist $x_0^{(k)}\in\F_{I}^{(n,k)}\ (0\leq k\leq n-1)$ 
and $x_0^{(n)}\in\G_{I}^{(n,n)}$ such that $x_0=\sum_{k=0}^nx_0^{(k)}$. 
We have $\|x_0\|=\max\{\|x_0^{(0)}\|,\ldots,\|x_0^{(n)}\|\}$.

First we consider the case that $\|x_0\|=\|x_0^{(k)}\|$ 
for some $k\leq n-1$.
If we set $q_k=\sum_{\mu\in\W_n^{(k)}}S_\mu p_nS_\mu^*\in M((\cpi)/I)$,
then $q_k$ is a projection satisfying that 
$q_kS_{\mu_l}S_{\nu_l}^*q_k=0$ if $|\mu_l|\neq |\nu_l|$.
Hence $q_kxq_k=q_kx_0q_k=x_0^{(k)}$.
We get $\|x_0\|=\|x_0^{(k)}\|=\|q_kxq_k\|\leq\|x\|$.
Next we consider the case that $\|x_0\|=\|x_0^{(n)}\|$.
Then there exists $\gamma_0\in X_I$ 
such that $\|x_0^{(n)}\|=\|x_0^{(n)}(\gamma_0)\|$.
By Lemma \ref{X^{(n)}} (iv), we have
$$X_I=\bigcup_{\mu\in\W_n}(X_I^{(n)}+\omega_\mu)\cup
\bigcap_{k=1}^\infty \bigg( \bigcup_{\mu\in\W_n^{(k)}}(X_I+\omega_\mu)\bigg).$$
When $\gamma_0\in X_I^{(n)}+\omega_\mu$ for some $\mu\in\W_n$, 
set $u=\sum_{\nu\in\W_n^{(n)}}S_\nu S_\mu p_nS_\nu^*\in M((\cpi)/I)$.
Then $u$ is a partial isometry.
We have $u^*xu=u^*x_0u=u^*x_0^{(n)}u=\pi_n(\sigma_{\omega_\mu}(x_0^{(n)}))$
where $\pi_n$ is the natural surjection 
from $\G_I^{(n,n)}$ onto $\F_I^{(n,n)}$.
Since $\gamma_0-\omega_\mu\in X_I^{(n)}$, we have 
$$\|\pi_n(\sigma_{\omega_\mu}(x_0^{(n)}))\|
\geq\|\sigma_{\omega_\mu}(x_0^{(n)})(\gamma_0-\omega_\mu)\|
=\|x_0^{(n)}(\gamma_0)\|=\|x_0^{(n)}\|=\|x_0\|.$$
Therefore $\|x_0\|\leq\|u^*xu\|\leq\|x\|$.

When $\gamma_0\in\bigcap_{k=1}^\infty 
\big( \bigcup_{\mu\in\W_n^{(k)}}(X_I+\omega_\mu)\big)$, 
we can find $i\in\{1,2,\ldots,n\}$ such that 
$\gamma_0-k\omega_i\in X_I$ for all $k\in\N$.
Since $\omega$ satisfies Condition \ref{cond}, 
either $k\omega_i\neq 0$ for any $k\in\Z_+$ or there exists a sequence 
$\{\mu_k\}_{k\in\Z_+}\subset\W_n$ with $\lim_{k\to\infty}\omega_{\mu_k}=0$ 
and $S_{\mu_k}^*S_i=0$ for any $k$.
In the case that $k\omega_i\neq 0$ for any $k\in\Z_+$, 
we can find a neighborhood $U$ of $\gamma_0-n\omega_i\in X_I$ such that 
$U\cap (U+k\omega_i)=\emptyset$ for $k=1,2,\ldots,n$.
Choose a function $f$ with $0\leq f\leq 1$ satisfying that 
the support of $f$ is contained in $U$ and $f(\gamma_0-n\omega_i)=1$.
Set $u=\sum_{\mu\in\W_n^{(n)}}S_\mu S_i^n f^{1/2} S_\mu^*\in(\cpi)/I$.
Since 
$$u^*u=\sum_{\mu,\nu\in\W_n^{(n)}}S_\mu f^{1/2} {S_i^*}^n S_\mu^* 
S_\nu S_i^n f^{1/2}S_\nu^*=\sum_{\mu\in\W_n^{(n)}}S_\mu f S_\mu^*,$$
$u^*u$ corresponds to $1\otimes f$ 
under the isomorphism $\G_I^{(n,n)}\cong\M_{n^n}\otimes C_0(X_I)$.
Thus we have $\|u^*u\|=\sup_{\gamma\in X_I}|f(\gamma)|=1$, 
and so $\|u\|=1$.
A routine computation shows that 
$u^*xu=u^*x_0^{(n)}u=f\sigma_{n\omega_i}x_0^{(n)}\in C_0(X_I,\M_{n^n})$.
Since $\gamma_0-n\omega_i\in X_I$, we have 
$$\|u^*xu\|\geq\|f(\gamma_0-n\omega_i)\sigma_{n\omega_i}x_0^{(n)}(\gamma_0-n\omega_i)\|=\|x_0^{(n)}(\gamma_0)\|=\|x_0\|.$$
Hence $\|x_0\|\leq\|u^*xu\|\leq\|x\|$.
Finally, we consider the case that there exists a sequence 
$\{\mu_k\}_{k\in\Z_+}\subset\W_n$ 
with $\lim_{k\to\infty}\omega_{\mu_k}=0$ 
and $S_{\mu_k}^*S_i=0$ for any $k\in\Z_+$.
For $k\in\Z_+$, define a partial isometry $u_k=\sum_{\mu\in\W_n^{(n)}}S_\mu S_i^nS_{\mu_k}S_\mu^*\in\Oi\subset M((\cpi)/I)$.
A routine computation shows that 
$u_k^*xu_k=u_k^*x_0^{(n)}u_k
=\sigma_{n\omega_i+\omega_{\mu_k}}x_0^{(n)}\in C_0(X_I,M_{n^n})$.
Since $\gamma_0-n\omega_i\in X_I$, we have
$$\|u_k^*xu_k\|\geq\|\sigma_{n\omega_i+\omega_{\mu_k}}x_0^{(n)}(\gamma_0-n\omega_i)\|=\|x_0^{(n)}(\gamma_0+\omega_{\mu_k})\|.$$
Hence we have 
$\|x_0^{(n)}(\gamma_0+\omega_{\mu_k})\|\leq\|u_k^*xu_k\|\leq\|x\|$ 
for any $k\in\Z_+$.
Therefore $\|x_0\|=\|x_0^{(n)}(\gamma_0)\|=\lim_{k\to\infty}\|x_0^{(n)}(\gamma_0+\omega_{\mu_k})\|\leq\|x\|$.

Hence the map 
\begin{align*}
\spa\{S_\mu fS_\nu^*&\mid\mu,\nu\in\W_\infty,\ f\in C_0(X_I)\}\ni x\\
&\mapsto\ x_0\in \spa\{S_\mu fS_\nu^*\mid\mu,\nu\in\W_\infty,\ |\mu|=|\nu|,\ f\in C_0(X_I)\}.
\end{align*}
is well-defined and norm-decreasing.
The extension $E_I$ of the map above is the desired conditional expectation 
onto $\F_I$.
Uniqueness is easy to verify.
\eprf

By uniqueness, the conditional expectation $E_{I}$ above coincides 
with the one in Lemma \ref{cond.exp1} when $I$ is gauge invariant.
Actually an ideal of $\cpi$ is gauge invariant 
if there exists such a conditional expectation,
as we see in the proof of the following theorem.

\bthm\label{idestr1}
Suppose that $\omega$ satisfies Condition \ref{cond}.
Then for any ideal $I$ we have $I_{\widetilde{X}_I}=I$, 
and so $I$ is gauge invariant.
Hence there is a one-to-one correspondence between the set of ideals of $\cpi$ 
and the set of $\omega$-invariant pairs of subsets of $\Gamma$.
\ethm

\bprf
If $X_I=\emptyset$, then $I=\cpi$ so $I_{\widetilde{X}_I}=I$.
Let $I$ be an ideal that is not $\cpi$, and set $J=I_{\widetilde{X}_I}$.
By the same way as in the proof of Proposition \ref{unique}, 
we can find a surjective $*$-homomorphism $\pi:(\cpi)/J\to(\cpi)/I$ 
whose restriction on $\F_{J}$ is an isomorphism from 
$\F_{J}$ onto $\F_{I}$.
By Proposition \ref{cond.exp2}, 
there exists a conditional expectation $E_I:(\cpi)/I\to\F_{I}$ 
satisfying $E_I\circ \pi=\pi\circ E_{J}$, 
where $E_J:(\cpi)/I\to\F_{J}$ is a faithful conditional expectation 
defined in Lemma \ref{cond.exp1}.
By Lemma \ref{isom}, $\pi$ is injective.
Therefore $I=I_{\widetilde{X}_I}$.
The last part follows from Theorem \ref{OneToOne}.
\eprf

\bco\label{primitive1}
When $\omega$ satisfies Condition \ref{cond}, 
an ideal $I$ of $\cpi$ is primitive 
if and only if the $\omega$-invariant pair $\widetilde{X}_I$ is prime.
\eco

\bprf
It follows from Proposition \ref{primepair} and Theorem \ref{idestr1}.
\eprf

\section{The ideal structure of $\cpi$ (part 2)}

In this section, we investigate the ideal structure of $\cpi$ 
when $\omega$ does not satisfy Condition \ref{cond}
i.e.\ there exists $i\in\Z_+$ such that 
$k\omega_i=0$ for some positive integer $k$, 
and that there exist no sequences $\mu_1,\mu_2,\ldots$ in $\W_\infty$ such that
$S_{\mu_k}^*S_i=0$ for any $k$ and $\lim_{k\to\infty}\omega_{\mu_k}=0$.
Note that such $i$ is unique.
Without loss of generality, we may assume $i=1$.
Let $K$ be the smallest positive integer satisfying $K\omega_1=0$.
Denote by $\Gamma'$ the quotient of $\Gamma$ 
by the subgroup generated by $\omega_1$, which is isomorphic to $\Z/K\Z$.
We denote by $[\gamma]$ and $[U]$ the images in $\Gamma'$ 
of $\gamma\in\Gamma$ and $U\subset\Gamma$ respectively.
We use the symbol $([\gamma],\theta)$ 
for denoting elements of $\Gamma'\times\T$.
Define $A=\cspa\{S_1^kf{S_1^*}^l\mid f\in C_0(\Gamma), k,l\in\N\}$ 
which is a $C^*$-subalgebra of $\cpi$.
In \cite{Ka1}, 
we defined a $C^*$-algebra $T_K$ and a continuous family of $*$-homomorphisms 
$\varphi_\gamma:A\to T_K$ for $\gamma\in\Gamma$.
Note that $\varphi_\gamma(x)=0$ 
if and only if $\varphi_{\gamma+\omega_1}(x)=0$ for $x\in A$.
We also defined $\psi_{\gamma,\theta}=\pi_\theta\circ \varphi_\gamma$
for $(\gamma,\theta)\in\Gamma\times\T$, 
where $\pi_\theta:T_K\to\M_{K}$ is a continuous family of $*$-homomorphisms.

\bde
For an ideal $I$ of $\cpi$, 
we define the closed subset $Y_I$ of $\Gamma'\times\T$ by
$$Y_I=\{([\gamma],\theta)\in\Gamma'\times\T\mid \psi_{\gamma,\theta}(x)=0
\mbox{ for all }x\in A\cap I\}.$$
We denote by $\widetilde{Y}_I$ the pair $(Y_I,X_I^{(\infty)})$ of 
a subset $Y_I$ of $\Gamma'\times\T$ 
and a subset $X_I^{(\infty)}$ of $\Gamma$.
\ede

\bde\label{omegaY}
For a pair $\widetilde{Y}=(Y,X^{(\infty)})$ of 
a subset $Y$ of $\Gamma'\times\T$ 
and a subset $X^{(\infty)}$ of $\Gamma$, 
we define subsets $X$ and $X^{(n)}$ of $\Gamma$ by 
\begin{align*}
X&=\{\gamma\in\Gamma\mid ([\gamma],\theta)\in Y
\mbox{ for some } \theta\in\T\},\\
X^{(n)}&=X^{(\infty)}\cup\bigcup_{i=n+1}^\infty(X+\omega_i).
\end{align*}

With this notation, 
a pair $\widetilde{Y}=(Y,X^{(\infty)})$ is called {\em $\omega$-invariant} 
if $(X,X^{(\infty)})$ is an $\omega$-invariant pair of subsets of $\Gamma$ 
and if $Y$ is a closed set satisfying that $[X^{(1)}]\times\T\subset Y$.
\ede

\bpr
For an ideal $I$ of $\cpi$, 
the pair $\widetilde{Y}_I$ is $\omega$-invariant.
\epr

\bprf
By \cite[Proposition 5.15]{Ka1}, we have
$$X_I=\{\gamma\in\Gamma\mid 
([\gamma],\theta)\in Y_I\mbox{ for some } \theta\in\T\}.$$
By the argument in the proof of \cite[Lemma 5.21]{Ka1}, we have
$$X_I^{(1)}=\{\gamma\in\Gamma\mid
\varphi_{\gamma}(x)=0\mbox{ for any } x\in A\cap I\}.$$
Therefore $[X_I^{(1)}]\times\T\subset Y_I$.
Thus the pair $\widetilde{Y}_I$ is $\omega$-invariant.
\eprf

We get the $\omega$-invariant pair $\widetilde{Y}_I$ 
from an ideal $I$ of $\cpi$.
Conversely, from an $\omega$-invariant pair $\widetilde{Y}$ , 
we can construct the ideal $I_{\widetilde{Y}}$ of $\cpi$.

\bde
For an $\omega$-invariant pair $\widetilde{Y}=(Y,X^{(\infty)})$, 
we define $J_{\widetilde{Y}}\subset A$ 
and $I_{\widetilde{Y}}\subset\cpi$ by
\begin{align*}
J_{\widetilde{Y}}&=\{x\in A\mid 
\psi_{\gamma,\theta}(x)=0\mbox{ for }([\gamma],\theta)\in Y,
\mbox{ and } \varphi_\gamma(x)=0\mbox{ for }\gamma\in X^{(1)}\},\\
I_{\widetilde{Y}}
&=\cspa\big(\{S_\mu xS_\nu^*\mid \mu,\nu\in\W_\infty,\ 
x\in J_{\widetilde{Y}}\}\\
&\hspace*{2cm}\cup\{S_\mu p_nfS_\nu \mid \mu,\nu\in\W_\infty,\ 
n\in\Z_+,\ f\in C_0(\Gamma\setminus X^{(n)})\}\big),
\end{align*}
with the notation in Definition \ref{omegaY}.
\ede

\bpr
For an $\omega$-invariant pair $\widetilde{Y}$, 
$I_{\widetilde{Y}}$ is an ideal of $\cpi$.
\epr

\bprf
Once noting that 
$J_{\widetilde{Y}}\cap C_0(\Gamma)=C_0(\Gamma\setminus X)$ and 
$J_{\widetilde{Y}}\cap p_1 C_0(\Gamma)=p_1 C_0(\Gamma\setminus X^{(1)})$,
we can prove that $I_{\widetilde{Y}}$ is an ideal 
in a similar way to Proposition \ref{I_X,Xinfty} 
with the help of the computation in \cite[Proposition 5.20]{Ka1}.
\eprf

\blem\label{Y_I_Y}
Let $\widetilde{Y}=(Y,X^{(\infty)})$ be an $\omega$-invariant pair. 
For any $([\gamma],\theta)\not\in Y$, there exists 
$x\in J_{\widetilde{Y}}$ such that $\psi_{\gamma,\theta}(x)\neq 0$.
\elem

\bprf
The proof goes exactly the same as in the proof of \cite[Lemma 5.22]{Ka1}, 
once noting that $([\gamma],\theta)\not\in Y$ implies $\gamma\not\in X^{(1)}$. 
\eprf

\bpr\label{exist2}
Let $\widetilde{Y}=(Y,X^{(\infty)})$ be an $\omega$-invariant pair, 
and set $I=I_{\widetilde{Y}}$.
Then we have $\widetilde{Y}_I=\widetilde{Y}$.
\epr

\bprf
By Lemma \ref{Y_I_Y}, we get $Y_I\subset Y$.
To prove the other inclusion, it is sufficient to see 
that $\psi_{\gamma,\theta}(x)=0$ 
for $([\gamma],\theta)\in Y$ and $x\in I\cap A$. 
Take $\e>0$ arbitrarily. 
Since $x\in I$,
there exist $\mu_l,\nu_l\in\W_\infty$, $x_l\in J_{\widetilde{Y}}$ 
for $l=1,2,\ldots,L$ and $\mu_k',\nu_k'\in\W_\infty$, $n_k\in\Z_+$, 
$f_k\in C_0(\Gamma\setminus X^{(n_k)})$ for $k=1,2,\ldots,K$ 
such that 
$$\bigg\|x-\sum_{l=1}^LS_{\mu_l}x_lS_{\nu_l}^*-\sum_{k=1}^K
S_{\mu_k'}p_{n_k}f_kS_{\nu_k'}^*\bigg\|<\e.$$
Take a positive integer $m$ such that $m\geq |\mu_l|,|\nu_l|$ 
for any $l$ and $m> |\mu_k'|,|\nu_k'|$ for any $k$.
Then, $\left\|{S_1^*}^mxS_1^m-\sum_{l=1}^Lx_l'\right\|<\e$ 
where $x_l'={S_1^*}^mS_{\mu_l}x_lS_{\nu_l}^*S_1^m$ for $l=1,2,\ldots,L$.
Since $x_l'\in J_{\widetilde{Y}}$, 
we have $\|\psi_{\gamma,\theta}({S_1^*}^mxS_1^m)\|<\e$.
Since $\psi_{\gamma,\theta}(S_1)$ is a unitary, 
we have $\|\psi_{\gamma,\theta}(x)\|<\e$ for arbitrary $\e>0$.
Hence, we have $\psi_{\gamma,\theta}(x)=0$.
Therefore we get $Y_I=Y$.

From $Y_I=Y$, we have $X_I=X$.
By the definition of $I$, 
we see that $X_I^{(n)}\subset X^{(n)}$ for $n\in\Z_+$.
To the contrary, assume that $X_I^{(n)}\subsetneqq X^{(n)}$.
Then there exists $f\in C_0(\Gamma)$ 
such that $p_nf\in I$ and $f(\gamma_0)=1$ for some $\gamma_0\in X^{(n)}$.
Since $p_nf\in I$, there exist 
$\mu_l,\nu_l\in\W_\infty$, $x_l\in J_{\widetilde{Y}}$ 
for $l=1,2,\ldots,L$ and $\mu_k',\nu_k'\in\W_\infty$, $n_k\in\Z_+$, 
$f_k\in C_0(\Gamma\setminus X^{(n_k)})$ for $k=1,2,\ldots,K$ 
such that 
$$\bigg\|p_nf-\sum_{l=1}^LS_{\mu_l}x_lS_{\nu_l}^*-\sum_{k=1}^K
S_{\mu_k'}p_{n_k}f_kS_{\nu_k'}^*\bigg\|<\frac12.$$
Take a positive integer $m$ so large that 
$\mu_l,\nu_l,\mu_k',\nu_k'\in\W_{m}$, $n_k\leq m$ 
for any $l,k$ and $n\leq m$.
By Lemma \ref{X^{(n)}} (iii), we have 
$X^{(n)}=X^{(m)}\cup\bigcup_{i=n+1}^{m}(X+\omega_i)$.
We first consider the case that $\gamma_0\in X^{(m)}$.
By \cite[Lemma 5.4]{Ka1}, there exists 
$g_l\in C_0(\Gamma\setminus X^{1})$ with $p_1x_lp_1=p_1g_l$ for any $l$.
Hence we have $p_{m}x_lp_{m}=p_{m}p_1x_lp_1p_{m}=p_{m}g_l$ for any $l$.
Since
$$p_{m}\bigg(p_nf-\sum_{l=1}^LS_{\mu_l}x_lS_{\nu_l}^*
-\sum_{k=1}^K S_{\mu_k'}p_{n_k}f_kS_{\nu_k'}^*\bigg)p_{m}=
p_{m}f-\sum_{\mu_l=\nu_l=\emptyset}p_{m}g_l
-\sum_{\mu_k'=\nu_k'=\emptyset}p_{m}f_k,$$
we get 
$\|f-\sum_{\mu_l=\nu_l=\emptyset}g_l-\sum_{\mu_k'=\nu_k'=\emptyset}f_k\|<1/2$.
This contradicts the fact that $f(\gamma_0)=1$, $g_l(\gamma_0)=0$ and 
$f_k(\gamma_0)=0$ for any $l,k$.
When $\gamma_0\in X+\omega_i$ for some $i$ with $n<i\leq m$, 
we have $\sigma_{\omega_i}f=S_i^*(p_nf)S_i\in I$ and 
$\sigma_{\omega_i}f(\gamma_0-\omega_i)=1$.
This contradicts the fact that $X_I=X$.
Therefore $X_I^{(n)}=X^{(n)}$ for a positive integer $n$.
Hence $X_I^{(\infty)}=\bigcap_{n=1}^\infty X_I^{(n)}=
\bigcap_{n=1}^\infty X^{(n)}=X^{(\infty)}$.
Thus we have $\widetilde{Y}_I=\widetilde{Y}$.
\eprf

\bco\label{Ysub}
For two $\omega$-invariant pairs $\widetilde{Y}_1=(Y_1,X_1^{(\infty)})$, 
$\widetilde{Y}_2=(Y_2,X_2^{(\infty)})$,
we have $I_{\widetilde{Y}_1}\subset I_{\widetilde{Y}_2}$
if and only if $Y_1\supset Y_2$ and $X_1^{(\infty)}\supset X_2^{(\infty)}$.
\eco

A relation between $I_{\widetilde{Y}}$ and $I_{\widetilde{X}}$ 
can be described as follows.

\bpr\label{rotation}
Let $\widetilde{Y}=(Y,X^{(\infty)})$ be an $\omega$-invariant pair.
For $t\in\T$, set $\widetilde{Y}_t=(Y_t,X^{(\infty)})$ where 
$Y_t=\{([\gamma],\theta)\in\Gamma'\times\T\mid ([\gamma],t\theta)\in Y\}$.
Then $\widetilde{Y}_t$ is $\omega$-invariant 
and $\beta_t(I_{\widetilde{Y}})=I_{\widetilde{Y}_{t^K}}$ 
where $\beta$ is the gauge action.
We also have $I_{\widetilde{X}}
=\bigcap_{t\in\T}I_{\widetilde{Y}_t}$
where $\widetilde{X}=(X,X^{(\infty)})$ and 
$X=\{\gamma\in\Gamma\mid ([\gamma],\theta)\in Y 
\mbox{ for some }\theta\in\T\}$.
\epr

\bprf
See \cite[Proposition 5.24]{Ka1}.
\eprf

\bpr
For an $\omega$-invariant pair $\widetilde{X}=(X,X^{(\infty)})$ 
of subsets of $\Gamma$,
the pair $\widetilde{Y}=([X]\times\T,X^{(\infty)})$ is $\omega$-invariant
and $I_{\widetilde{Y}}=I_{\widetilde{X}}$.
\epr

\bprf
Obvious by Proposition \ref{rotation}.
\eprf

Now, we turn to showing that $I_{\widetilde{Y}_I}=I$ for any ideal $I$ 
(Theorem \ref{idestr2}).
To see this, we examine the primitive ideal space of $\cpi$.
Set $\vcsg=\csg\setminus\{0,\omega_1,\ldots,(K-1)\omega_1\}$.

\blem\label{isolate}
We have $\vcsg=\overline{\bigcup_{i=2}^\infty(\sg+\omega_i)}$ and 
$\vcsg$ is an $\omega$-invariant set.
\elem

\bprf
For $\gamma\in\csg$, we can find $\mu_k\in\W_\infty$ such that 
$\gamma=\lim_{k\to\infty}\omega_{\mu_k}$.
If $\mu_k=(1,1,\ldots,1)$ for sufficiently large $k$,
then $\gamma=m\omega_1$ for some $m\in\N$.
Hence for $\gamma\in\vcsg$, we can find $\mu_k\in\W_\infty$ 
with $\omega_{\mu_k}\in\bigcup_{i=2}^\infty(\sg+\omega_i)$ such that 
$\gamma=\lim_{k\to\infty}\omega_{\mu_k}$.
Thus $\vcsg\subset\overline{\bigcup_{i=2}^\infty(\sg+\omega_i)}$.
To prove the other inclusion, 
suppose $m\omega_1\in\overline{\bigcup_{i=2}^\infty(\sg+\omega_i)}$ 
for some $0\leq m<K$
and we will derive a contradiction.
In this case, $0$ is also in $\overline{\bigcup_{i=2}^\infty(\sg+\omega_i)}$.
Hence there exists a sequence $\{\mu_k\}$ in $\W_\infty$ 
with $S_{\mu_k}^*S_1=0$ 
such that $0=\lim_{k\to\infty}\omega_{\mu_k}$.
This contradicts the fact that $\omega$ does not satisfy Condition \ref{cond}.
Therefore $\vcsg=\overline{\bigcup_{i=2}^\infty(\sg+\omega_i)}$. 
From this equality, it is easy to see that 
$\vcsg$ is an $\omega$-invariant set.
\eprf

\bco\label{cptnbhd}
For any $\gamma_0\in\Gamma$, 
there exists a compact neighborhood $X$ of $\gamma_0$ satisfying that
$X\cap(X+\gamma)=\emptyset$ for any $\gamma\in\csg\setminus\{0\}$.
\eco

\bprf
Since 
$\csg\setminus\{0\}=\vcsg\cup\{\omega_1,2\omega_1,\ldots,(K-1)\omega_1\}$ 
is closed by Lemma \ref{isolate}, there exists a neighborhood $U$ 
of $0$ with $U\cap (\csg\setminus\{0\})=\emptyset$.
If we take a compact neighborhood $V$ of $0$ such that $V-V\subset U$,
then $X=\gamma_0+V$ becomes a desired compact neighborhood of $\gamma_0$.
\eprf

\blem\label{H_X}
For an $\omega$-invariant set $X$, we have 
$H_X=\bigcap_{n=1}^\infty\overline{\bigcup_{i=n}^\infty(X+\omega_i)}$.
If two $\omega$-invariant sets $X_1$ and $X_2$ satisfy $X_1\subset X_2$, 
then $H_{X_1}\subset H_{X_2}$.
\elem

\bprf
The former part follows from $X=X+\omega_1$, and this implies the latter part. 
\eprf

\bpr
For any $\gamma\in\Gamma$, we have $\gamma\notin H_{\gamma+\csg}$.
\epr

\bprf
By Lemma \ref{H_X}, we have 
$$H_{\gamma+\csg}=
\bigcap_{n=1}^\infty\overline{\bigcup_{i=n}^\infty(\gamma+\csg+\omega_i)}
\subset\overline{\bigcup_{i=2}^\infty(\gamma+\csg+\omega_i)}
=\gamma+\vcsg.$$
Hence $\gamma\notin H_{\gamma+\csg}$.
\eprf

For $\gamma\in\Gamma$, we set $P_{\gamma}=I_{\widetilde{X}}$ where 
$\widetilde{X}=(\gamma+\csg,H_{\gamma+\csg}\cup\{\gamma\})$ 
which is a prime $\omega$-invariant pair.
We will show that $P_{\gamma}$ is the unique primitive ideal satisfying that 
$\widetilde{X}_{P_{\gamma}}=(\gamma+\csg,H_{\gamma+\csg}\cup\{\gamma\})$.
To see this, we need the following lemma.

\blem\label{Pgamma0}
Let $I$ be an ideal of $\cpi$ with $X_I=\overline{X_I^{(\infty)}+\sg}$.
Then $I=I_{\widetilde{X}_I}$.
\elem

\bprf
By the argument in the proof of 
Proposition \ref{cond.exp2} and Theorem \ref{idestr1},
it suffices to show that $\|x_0\|\leq\|x\|$ for 
$x=\sum_{l=1}^LS_{\mu_l}f_lS_{\nu_l}^*\in (\cpi)/I$ and 
$x_0=\sum_{|\mu_l|=|\nu_l|}S_{\mu_l}f_lS_{\nu_l}^*$.
If we choose a positive integer $n$ 
so that $|\mu_l|,|\nu_l|\leq n$ and $\mu_l,\nu_l\in\W_n$ for any $l$, 
then $x_0\in \F_I^{(n)}$.
We can find $x_0^{(k)}\in\F_{I}^{(n,k)}\ (0\leq k\leq n-1)$ 
and $x_0^{(n)}\in\G_{I}^{(n,n)}$ such that $x_0=\sum_{k=0}^nx_0^{(k)}$. 
We have $\|x_0\|=\max\{\|x_0^{(0)}\|,\ldots,\|x_0^{(n)}\|\}$.
In the case that $\|x_0\|=\|x_0^{(k)}\|$ for some $k\leq n-1$, 
we can prove $\|x_0\|\leq\|x\|$ in a similar way to the proof of 
Proposition \ref{cond.exp2}.
In the case that $\|x_0\|=\|x_0^{(n)}\|$, 
there exists $\gamma_0\in X_I$ 
such that $\|x_0\|=\|x_0^{(n)}(\gamma_0)\|$.
Since $X_I=\overline{X_I^{(\infty)}+\sg}$,
there exist a sequence $\mu_1,\mu_2,\ldots\in\W_\infty$ 
and a sequence $\gamma_1,\gamma_2,\ldots,\in X_I^{(\infty)}$
such that $\gamma_0=\lim_{k\to\infty}(\gamma_k+\omega_{\mu_k})$.
We can find sequences $\mu_1',\mu_2',\ldots\in\W_\infty$ and 
$\nu_1,\nu_2,\ldots\in\W_n$ such that 
$\omega_{\mu_k}=\omega_{\mu_k'}+\omega_{\nu_k}$ and
none of $1,2,\ldots,n$ appears in the word $\mu_k'$ for any $k$.
For $k\in\Z_+$, define a partial isometry 
$u_k=\sum_{\mu\in\W_n^{(n)}}S_\mu S_{\nu_k}p_n S_\mu^*$.
We have 
$u_k^*xu_k=u_k^*x_0u_k=u_k^*x_0^{(n)}u_k
=\pi_n(\sigma_{\omega_{\nu_k}}x_0^{(n)})$, 
where $\pi_n$ is the natural surjection 
from $\G_I^{(n,n)}$ onto $\F_I^{(n,n)}$.
Since $\gamma_k\in X_I^{(\infty)}$, 
we have $\gamma_k+\omega_{\mu_k'}\in X_I^{(n)}$.
Hence 
$$\|\pi_n(\sigma_{\omega_{\nu_k}}x_0^{(n)})\|
\geq\|\sigma_{\omega_{\nu_k}}x_0^{(n)}(\gamma_k+\omega_{\mu_k'})\|
=\|x_0^{(n)}(\gamma_k+\omega_{\mu_k'}+\omega_{\nu_k})\|.$$
Therefore we get 
$$\|x_0\|=\|x_0^{(n)}(\gamma_0)\|=\lim_{k\to\infty}
\|x_0^{(n)}(\gamma_k+\omega_{\mu_k'}+\omega_{\nu_k})\|\leq\|x\|.$$
We are done.
\eprf

\bpr\label{Pgamma}
For any $\gamma\in\Gamma$, the ideal $P_{\gamma}$ is the unique primitive 
ideal satisfying that 
$\widetilde{X}_{P_{\gamma}}=(\gamma+\csg,H_{\gamma+\csg}\cup\{\gamma\})$.
\epr

\bprf
To prove that $P_{\gamma}$ is primitive, 
it suffices to show that it is prime because $\cpi$ is separable.
Let $I_1,I_2$ be ideals of $\cpi$ 
with $I_1\cap I_2= P_\gamma$.
Then we get 
$\widetilde{X}_{I_1}\cup\widetilde{X}_{I_2}=\widetilde{X}_{P_\gamma}$.
Since $\widetilde{X}_{P_\gamma}$ is a prime $\omega$-invariant pair, 
we have either 
$\widetilde{X}_{I_1}=\widetilde{X}_{P_\gamma}$ or 
$\widetilde{X}_{I_2}=\widetilde{X}_{P_\gamma}$.
By Lemma \ref{Pgamma0}, we have either $I_1=P_\gamma$ or $I_2=P_\gamma$.
Therefore $P_\gamma$ is primitive.
The uniqueness follows from Lemma \ref{Pgamma0}.
\eprf

We denote by $\varDelta$ the set of prime $\omega$-invariant sets 
which are not of the form $\gamma+\csg$.
For $X\in\varDelta$, we denote by $P_X$ the ideal $I_{\widetilde{X}}$
for $\widetilde{X}=(X,H_X)$ which is a prime $\omega$-invariant pair.
We will show that for any $X\in\varDelta$, 
$P_X$ is the unique primitive ideal satisfying $\widetilde{X}_{P_X}=(X,H_X)$.

\blem\label{PX0}
Let $X\in\varDelta$ and $\gamma\in X$.
Then there exist a sequence $\mu_1,\mu_2,\ldots$ in $\W_\infty$ 
and a sequence $\gamma_1,\gamma_2,\ldots$ in $X$ such that 
$S_{\mu_k}^*S_1=0$ for any $k$ 
and $\gamma=\lim_{k\to\infty}(\gamma_k+\omega_{\mu_k})$.
\elem

\bprf
Since $X\in\varDelta$, there exists $\gamma'\in X\setminus (\gamma+\csg)$.
Since $X$ is prime, Proposition \ref{Xprime} gives us 
two sequences $\mu_1,\mu_2,\ldots$, $\nu_1,\nu_2,\ldots$ in $\W_\infty$ 
and a sequence $\gamma_1,\gamma_2,\ldots$ in $X$ with 
$\gamma=\lim_{k\to\infty}(\gamma_k+\omega_{\mu_k})$
and $\gamma'=\lim_{k\to\infty}(\gamma_k+\omega_{\nu_k})$.
We will show that we can choose such $\mu_k$ satisfying $S_{\mu_k}^*S_1=0$.
If not so, then $\mu_k=(1,1,\ldots,1)$ for sufficiently large $k$.
This implies $\gamma'=\lim_{k\to\infty}(\gamma-|\mu_k|\omega_1+\omega_{\nu_k})$
which contradicts the fact that $\gamma'\notin\gamma+\csg$.
Therefore we can find desired sequences.
\eprf

\blem\label{PX1}
If an ideal $I$ of $\cpi$ satisfies $X_I\in\varDelta$,
then $I=I_{\widetilde{X}_I}$.
\elem

\bprf
Similarly as the proof of Lemma \ref{Pgamma0}, 
it suffices to show that $\|x_0\|\leq\|x\|$ for 
$x=\sum_{l=1}^LS_{\mu_l}f_lS_{\nu_l}^*\in (\cpi)/I$ and 
$x_0=\sum_{|\mu_l|=|\nu_l|}S_{\mu_l}f_lS_{\nu_l}^*\in \F_I^{(n)}$.
We can find $x_0^{(k)}\in\F_{I}^{(n,k)}\ (0\leq k\leq n-1)$ 
and $x_0^{(n)}\in\G_{I}^{(n,n)}$ such that $x_0=\sum_{k=0}^nx_0^{(k)}$. 
We have $\|x_0\|=\max\{\|x_0^{(0)}\|,\ldots,\|x_0^{(n)}\|\}$.
In the case that $\|x_0\|=\|x_0^{(k)}\|$ for some $k\leq n-1$, 
we can prove $\|x_0\|\leq\|x\|$ in a similar way to the proof of 
Proposition \ref{cond.exp2}.
In the case that $\|x_0\|=\|x_0^{(n)}\|$, 
there exists $\gamma_0\in X_I$ 
such that $\|x_0\|=\|x_0^{(n)}(\gamma_0)\|$.
By Lemma \ref{PX0}, we have a sequence $\mu_1,\mu_2,\ldots$ in $\W_\infty$ 
and a sequence $\gamma_1,\gamma_2,\ldots$ in $X_I$ such that 
$S_{\mu_k}^*S_1=0$ and $\gamma=\lim_{k\to\infty}(\gamma_k+\omega_{\mu_k})$.
For $k\in\Z_+$, set a partial isometry 
$u_k=\sum_{\mu\in\W_n^{(n)}}S_\mu S_1^{Kn}S_{\mu_k}S_\mu^*$.
We have 
$u_k^*xu_k=u_k^*x_0u_k=u_k^*x_0^{(n)}u_k=\sigma_{\omega_{\mu_k}}x_0^{(n)}$.
Since $\gamma_k\in X_I$, we have 
$$\|u_k^*xu_k\|\geq\|\sigma_{\omega_{\mu_k}}x_0^{(n)}(\gamma_k)\|
=\|x_0^{(n)}(\gamma_k+\omega_{\mu_k})\|.$$
Therefore we get 
$$\|x_0\|=\|x_0^{(n)}(\gamma)\|=\lim_{k\to\infty}
\|x_0^{(n)}(\gamma_k+\omega_{\mu_k})\|\leq\|x\|.$$
We are done.
\eprf

\bpr\label{PX}
For $X\in\varDelta$, the ideal $P_{X}$ is the unique primitive ideal 
satisfying $\widetilde{X}_{P_X}=(X,H_X)$.
\epr

\bprf
With the help of Lemma \ref{PX1},
the proof goes similarly as the one in Proposition \ref{Pgamma}.
\eprf

By Proposition \ref{pp}, the remaining candidates for primitive ideals
are ideals $P$ satisfying 
$\widetilde{X}_P=(\gamma_0+\csg,H_{\gamma_0+\csg})$
for some $\gamma_0\in\Gamma$.
We will determine such primitive ideals.

\bde\label{DefP}
For $([\gamma],\theta)\in\Gamma'\times\T$, 
we set $Y_{([\gamma],\theta)}=\big\{([\gamma],\theta)\big\}\cup
\big([\gamma+\vcsg]\times\T\big)$. 
Then $\widetilde{Y}=(Y_{([\gamma],\theta)},H_{\gamma+\csg})$ 
is an $\omega$-invariant pair.
We write $P_{([\gamma],\theta)}$ 
for denoting $I_{\widetilde{Y}}$.
\ede

We can show that $P_{([\gamma],\theta)}$ is a primitive ideal 
for any $([\gamma],\theta)\in\Gamma'\times\T$ 
by using the technique in \cite{Ka1}.
To do so, we need Proposition \ref{local}, which 
will also be used to determine the topology of primitive ideal space of $\cpi$.

\blem\label{I_X=}
For an $\omega$-invariant set $X$, 
the pair $\widetilde{X}=(X,X)$ is $\omega$-invariant and we have
$$I_{\widetilde{X}}=\cspa\{S_\mu fS_\nu^*\mid 
\mu,\nu\in\W_\infty,\ f\in C_0(\Gamma\setminus X)\}.$$
\elem

\bprf
Clearly, $\widetilde{X}=(X,X)$ is $\omega$-invariant. 
Set $I=\cspa\{S_\mu fS_\nu^*\mid 
\mu,\nu\in\W_\infty,\ f\in C_0(\Gamma\setminus X)\}$. 
In a similar way to the proof of Proposition \ref{I_X,Xinfty}, 
we can see that $I$ is a gauge-invariant ideal of $\cpi$. 
We also see that $X_I^{(n)}=X$ for any $n\in\N$
by arguing as in the proof of Proposition \ref{exist}.
Hence $I_{\widetilde{X}}=I$ by Theorem \ref{OneToOne}.
\eprf

\bpr\label{local}
Let $X$ be a compact subset of $\Gamma$ such that 
$X\cap (X+\gamma)=\emptyset$ for any $\gamma\in\csg\setminus\{0\}$,
and set $X_1=X+\csg$ and $X_2=X+\vcsg$. 
Then we have that $\widetilde{X}_0=(X_1,X_1)$, $\widetilde{X}_1=(X_1,X_2)$ and 
$\widetilde{X}_2=(X_2,X_2)$ are $\omega$-invariant pairs,
and that
\begin{align*}
I_{\widetilde{X}_2}/I_{\widetilde{X}_1}&\cong\K\otimes C(X\times\T),&
I_{\widetilde{X}_1}/I_{\widetilde{X}_0}&\cong\K\otimes C(X_1\setminus X_2).
\end{align*}
\epr

\bprf
Since $X$ is compact and $\csg$ is closed, $X_1=X+\csg$ becomes closed. 
The same reason shows that $X_2$ is closed.
By Lemma \ref{isolate}, both $X_1$ and $X_2$ are $\omega$-invariant and
$X_2=\overline{\bigcup_{i=2}^\infty(X_1+\omega_i)}$.
Therefore $\widetilde{X}_0,\widetilde{X}_1,\widetilde{X}_2$ are
$\omega$-invariant pairs.
Since $I_{\widetilde{X}_1}\cap p_1 C_0(\Gamma)=p_1 C_0(\Gamma\setminus X_2)$,
we have $p_1f=0$ 
for any $f\in C_0(X_1\setminus X_2)\subset I_{\widetilde{X}_2}/I_{\widetilde{X}_1}$. 
Note that $X_1\setminus X_2$ is a disjoint union of compact sets 
$X,X+\omega_1,\ldots,X+(K-1)\omega_1$. 
For $f\in C(X+m\omega_1)\subset I_{\widetilde{X}_2}/I_{\widetilde{X}_1}$ 
with $0<m<K$, 
we have $\sigma_{m\omega_1}f\in C(X)$ and 
\begin{align*}
S_1^m\sigma_{m\omega_1}f{S_1^*}^m
&=S_1^{m-1}S_1S_1^*\sigma_{(m-1)\omega_1}f{S_1^*}^{m-1}
=S_1^{m-1}\sigma_{(m-1)\omega_1}f{S_1^*}^{m-1}\\
&=\cdots=f.
\end{align*}
Hence, we have
$I_{\widetilde{X}_2}/I_{\widetilde{X}_1}
=\cspa\{S_\mu fS_\nu^*\mid \mu,\nu\in\W_\infty, f\in C(X)\}$ 
by Lemma \ref{I_X=}.
Set $\W_\infty^+=\W_\infty\setminus\{\mu 1^K\in\W_\infty\mid\mu\in\W_\infty\}$
and denote by $\chi$ the characteristic function of $X$.
Then $\{S_\mu\chi S_\nu^*\}_{\mu,\nu\in \W_\infty^+}$ 
satisfies the relation of matrix units 
and $\sum_{\mu\in\W_\infty^+}S_\mu\chi S_\mu^*=1$ 
(strictly).
Hence we have $I_{\widetilde{X}_2}/I_{\widetilde{X}_1}\cong\K\otimes B$ 
where $B=\chi(I_{\widetilde{X}_2}/I_{\widetilde{X}_1})\chi$.
We have 
$$B=\cspa\{\chi S_\mu fS_\nu^* \chi\mid \mu,\nu\in\W_\infty, f\in C(X)\}
=\cspa\{(S_1^K)^m f\mid m\in\Z, f\in C(X)\}.$$
Since $B$ is generated by $C(X)$ and a unitary $S_1^K\chi$ which commute 
with each other and since $B$ is globally invariant under the gauge action, 
we have $B\cong C(X\times\T)$.
Therefore we get 
$I_{\widetilde{X}_2}/I_{\widetilde{X}_1}\cong\K\otimes C(X\times\T)$.

By the definition, 
$$I_{\widetilde{X}_1}/I_{\widetilde{X}_0}
=\cspa\{S_\mu p_nfS_\nu^*\mid \mu,\nu\in\W_\infty, n\geq 1,\ 
f\in C(X_1\setminus X_2)\}.$$
For $f\in C(X_1\setminus X_2)\subset 
I_{\widetilde{X}_1}/I_{\widetilde{X}_0}$ and $i\geq 2$, 
we have $S_iS_i^*f=S_i\sigma_{\omega_i}fS_i^*=0$.
Hence $p_nf=p_1f$ for any $n\geq 1$ 
and any $f\in C(X_1\setminus X_2)$.
Thus $I_{X_1,X_2}/I_{X_1,X_1}
=\cspa\{S_\mu p_2fS_\nu^*\mid \mu,\nu\in\W_\infty,\ 
f\in C(X_1\setminus X_2)\}.$
We can show that 
$\{S_\mu p_2\chi'S_\nu^*\}_{\mu,\nu\in\W_\infty}$ 
satisfies the relation of matrix units 
and $\sum_{\mu\in\W_\infty}S_\mu p_2\chi'S_\mu^*=1$ (strictly),
where $\chi'$ is the characteristic function of $X_1\setminus X_2$.
Hence we have $I_{\widetilde{X}_1}/I_{\widetilde{X}_0}\cong\K\otimes B'$ 
where 
$$B'=p_2\chi'(I_{\widetilde{X}_1}/I_{\widetilde{X}_0})p_2\chi'
=\cspa\{p_2f\mid f\in C(X_1\setminus X_2)\}\cong C(X_1\setminus X_2).$$
Therefore we get 
$I_{\widetilde{X}_1}/I_{\widetilde{X}_0}\cong\K\otimes C(X_1\setminus X_2)$.
\eprf

With the help of Proposition \ref{local},
we have the following proposition 
by exactly the same argument as the proof of \cite[Proposition 5.41]{Ka1}.

\bpr\label{Primsurj}
For $\gamma_0\in\Gamma$, the set of all primitive ideals $P$ satisfying 
$\widetilde{X}_P=(\gamma_0+\csg,H_{\gamma_0+\csg})$ 
is $\{P_{([\gamma_0],\theta)}\mid\theta\in\T\}$.
\epr

Now, we can describe the primitive ideal space $\Prim(\cpi)$ of $\cpi$
as follows.

\bpr\label{primitive2}
We have 
$\Prim(\cpi)=\{P_z\mid z\in (\Gamma'\times\T)\sqcup\Gamma\sqcup\varDelta\}$,
where $\sqcup$ means a disjoint union.
\epr

The primitive ideal space $\Prim(\cpi)$ is a topological space 
whose closed sets are given by 
$\{P\in\Prim(\cpi)\mid I\subset P\}$ for ideals $I$.
We will investigate which subset of 
$(\Gamma'\times\T)\sqcup\Gamma\sqcup\varDelta$ 
corresponds to a closed subset of $\Prim(\cpi)$.
By Corollary \ref{Ysub}, the following is easy to verify.

\blem\label{Primlem1}
Let $\widetilde{Y}=(Y,X^{(\infty)})$ be an $\omega$-invariant set.
\benu
\item For $([\gamma],\theta)\in\Gamma'\times\T$, we have
$I_{\widetilde{Y}}\subset P_{([\gamma],\theta)}$ if and only if 
$([\gamma],\theta)\in Y$.
\item For $\gamma\in\Gamma$, we have 
$I_{\widetilde{Y}}\subset P_{\gamma}$ if and only if 
$\gamma\in X^{(\infty)}$.
\item For $X\in\varDelta$, we have 
$I_{\widetilde{Y}}\subset P_{X}$ if and only if 
$[X]\times\T\subset Y$.
\eenu
\elem

\blem\label{Primlem2}
Let $X$ be a compact subset of $\Gamma$ 
such that $X\cap (X+\gamma)=\emptyset$ for any $\gamma\in\csg\setminus\{0\}$,
and set $X_1=X+\csg$ and $X_2=X+\vcsg$, which are $\omega$-invariant sets.
If $X_0\in\varDelta$ satisfies $X_1\supset X_0$, then $X_2\supset X_0$.
\elem

\bprf
To the contrary, assume $X_0\in\varDelta$ satisfies 
$X_1\supset X_0$ and $X_2\not\supset X_0$.
Then $X_0\cap X\neq\emptyset$ and 
$(X_0\cap X)+\csg$ is an $\omega$-invariant set 
satisfying $(X_0\cap X)+\csg\subset X_0$.
Since $((X_0\cap X)+\csg) \cup X_2\supset X_0$ and $X_0$ is prime,
we have $(X_0\cap X)+\csg\supset X_0$.
Hence $X_0=(X_0\cap X)+\csg$.
If $X_0\cap X$ has two points $\gamma_1,\gamma_2$, 
then we can take open sets $U_1,U_2$ such that $\gamma_1\in U_1$, 
$\gamma_2\in U_2$, $U_1\cap U_2=\emptyset$.
Two $\omega$-invariant sets $X_1'=(X_0\cap X\setminus U_1)+\csg$, 
$X_2'=(X_0\cap X\setminus U_2)+\csg$ satisfies 
$X_1'\not\supset X_0$, $X_2'\not\supset X_0$ and $X_1'\cup X_2'=X_0$. 
This contradicts the primeness of $X_0$.
Hence $X_0\cap X$ is just a point.
However, this contradicts the fact that $X_0\in\varDelta$.
Therefore $X_2\supset X_0$ when $X_0\in\varDelta$ satisfies $X_1\supset X_0$.
\eprf

\blem\label{Primlem3}
Let $\widetilde{Y}_\lambda=(Y_\lambda,X_\lambda^{(\infty)})$ 
be an $\omega$-invariant pair for each $\lambda\in\Lambda$. 
Set $I=\bigcap_{\lambda\in\Lambda}I_{\widetilde{Y}_\lambda}$.
Then $Y_I=\overline{\bigcup_{\lambda\in\Lambda}Y_\lambda}$.
\elem

\bprf
For any $\lambda\in\Lambda$, we have $Y_I\supset Y_\lambda$
because $I\subset I_{\widetilde{Y}_\lambda}$.
Hence we get $Y_I\supset \overline{\bigcup_{\lambda\in\Lambda}Y_\lambda}$.
Take 
$([\gamma_0],\theta_0)\notin
\overline{\bigcup_{\lambda\in\Lambda}Y_\lambda}$.
Then 
there exists a neighborhood $U$ of $([\gamma_0],\theta_0)$ satisfying
$U\cap\overline{\bigcup_{\lambda\in\Lambda}Y_\lambda}=\emptyset$.
By the same argument as in the proof of \cite[Lemma 5.22]{Ka1},
we can find $x_0\in A$ such that 
$\psi_{([\gamma_0],\theta_0)}(x_0)\neq 0$ and 
$\psi_{([\gamma],\theta)}(x_0)=0$ if $([\gamma],\theta)\notin U$ 
and $\varphi_{\gamma}(x_0)=0$ if $([\gamma]\times\T)\cap U=\emptyset$.
Therefore we have $x_0\in I$, 
and it implies that 
$([\gamma_0],\theta_0)\notin Y_I$.
Thus $Y_I=\overline{\bigcup_{\lambda\in\Lambda}Y_\lambda}$.
\eprf

\blem\label{Primlem4}
For any $X\in\varDelta$, we have $P_{X}=\bigcap_{z\in [X]\times\T}P_z$.
\elem

\bprf
By Lemma \ref{Primlem1}, 
we have $P_{X}\subset\bigcap_{z\in [X]\times\T}P_z$.
By Lemma \ref{Primlem3}, 
we have $Y_{\bigcap_{z\in [X]\times\T}P_z}=[X]\times\T$. 
Hence we have $\bigcap_{z\in [X]\times\T}P_z\subset P_{X}$ by 
Lemma \ref{Primlem1}.
Thus $P_{X}=\bigcap_{z\in [X]\times\T}P_z$.
\eprf

In the proof of the following proposition,
we use the fact that 
the subset $\{P\in\Prim(\cpi)\mid I_1\subset P, I_2\not\subset P\}$ 
of $\Prim(\cpi)$ is homeomorphic to $\Prim(I_2/I_1)$,
for two ideals $I_1,I_2$ of $\cpi$ with $I_1\subset I_2$.

\bpr\label{closed}
Let $Z=Y\sqcup X^{(\infty)}\sqcup \varLambda$ be a subset of 
$(\Gamma'\times\T)\sqcup\Gamma\sqcup\varDelta$.
The set $P_Z=\{P_z\mid z\in Z\}$ is closed in $\Prim(\cpi)$ 
if and only if $(Y,X^{(\infty)})$ is an $\omega$-invariant set and
$\varLambda=\{X\in\varDelta\mid [X]\times\T\subset Y\}$.
\epr

\bprf
Let us take a subset $Z=Y\sqcup X^{(\infty)}\sqcup \varLambda$ of 
$(\Gamma'\times\T)\sqcup\Gamma\sqcup\varDelta$
satisfying that $(Y,X^{(\infty)})$ is an $\omega$-invariant set and
$\varLambda=\{X\in\varDelta\mid [X]\times\T\subset Y\}$.
Then the set $P_Z=\{P_z\mid z\in Z\}$ coincides 
with the closed subset defined by the ideal $I_{\widetilde{Y}}$
by Lemma \ref{Primlem1}. 

Conversely, assume $P_Z$ is closed, that is, 
there exists an ideal $I$ of $\cpi$ with 
$Z=\{z\in Y\sqcup X^{(\infty)}\sqcup \varLambda\mid I\subset P_z\}$.
We first show that $Y$ and $X^{(\infty)}$ is closed.
Take $\gamma_0\in\Gamma$ arbitrarily.
By Corollary \ref{cptnbhd}, 
there exists a compact neighborhood $X$ of $\gamma_0$ 
such that $X\cap (X+\gamma)=\emptyset$ for any $\gamma\in\csg\setminus\{0\}$.
Set $\widetilde{X}_0=(X_1,X_1)$, $\widetilde{X}_1=(X_1,X_2)$ and 
$\widetilde{X}_2=(X_2,X_2)$ where $X_1=X+\csg$ and $X_2=X+\vcsg$.
Note that $X\ni\gamma\mapsto[\gamma]\in [X_1\setminus X_2]$ is a homeomorphism.
By Lemma \ref{Primlem1} and Lemma \ref{Primlem2}, we have 
\begin{align*}
\big\{z\in (\Gamma'\times\T)\sqcup\Gamma\sqcup\varDelta\mid 
I_{\widetilde{X}_1}\subset P_z, I_{\widetilde{X}_2}\not\subset P_z\big\}
&=[X_1\setminus X_2]\times\T\subset \Gamma'\times\T,\\
\big\{z\in (\Gamma'\times\T)\sqcup\Gamma\sqcup\varDelta\mid 
I_{\widetilde{X}_0}\subset P_z, I_{\widetilde{X}_1}\not\subset P_z\big\}
&=X_1\setminus X_2\subset\Gamma.
\end{align*}
By Proposition \ref{local}, the map 
$[X_1\setminus X_2]\times\T\ni z\mapsto P_z\in\Prim(\cpi)$ is a homeomorphism 
from $[X_1\setminus X_2]\times\T$, 
whose topology is the relative topology of $\Gamma'\times\T$, 
to the subset 
$\{P\in\Prim(\cpi)\mid I_{\widetilde{X}_1}\subset P, 
I_{\widetilde{X}_2}\not\subset P\}$ 
of $\Prim(\cpi)$.
The set $Y\cap ([X_1\setminus X_2]\times\T)\subset \Gamma'\times\T$ is closed 
in $[X_1\setminus X_2]\times\T$ because $P_Y$ is closed.
Hence, the subset $Y$ is closed in $\Gamma'\times\T$.
Similarly $X^{(\infty)}$ is closed in $\Gamma$.
Set $X=\{\gamma\in\Gamma\mid ([\gamma],\theta)\in Y \mbox{ for some }\theta\in\T\}$, which is closed because $Y$ is closed.
Set $J=\bigcap_{([\gamma],\theta)\in Y}P_{([\gamma],\theta)}$.
We have $I\subset J$.
By Lemma \ref{Primlem3}, we have $Y_J=Y$.
Hence $H_X\subset X_J^{(\infty)}$.
We have $J\subset P_\gamma$ for any $\gamma\in H_X$ by Lemma \ref{Primlem1}.
Therefore $H_X\subset X^{(\infty)}$.
We have $([\gamma+\omega_i],\theta')\in Y$ 
for any $([\gamma],\theta)\in Y$, any $i\geq 2$ and any $\theta'\in\T$ 
because $P_{([\gamma],\theta)}\subset P_{([\gamma+\omega_i],\theta')}$.
Hence we get $[X+\omega_i]\times\T\subset Y$.
We also have $[X^{(\infty)}]\times\T\subset Y$ 
because $P_\gamma\subset P_{([\gamma],\theta)}$ 
for any $([\gamma],\theta)\in\Gamma'\times\T$.
Therefore we have proved that $(Y,X^{(\infty)})$ is an $\omega$-invariant set.
Finally, we have $\varLambda=\{X\in\varDelta\mid [X]\times\T\subset Y\}$
by Lemma \ref{Primlem4}.
It completes the proof.
\eprf

By the proposition above, we get the following.

\bthm\label{idestr2}
When $\omega$ does not satisfy Condition \ref{cond},
there is a one-to-one correspondence between the set of ideals of $\cpi$ 
and the set of $\omega$-invariant pairs of subsets of $\Gamma'\times\T$
and subsets of $\Gamma$.
Hence for any ideal $I$ of $\cpi$, we have $I=I_{\widetilde{Y}_I}$. 
\ethm

\bprf
There is a one-to-one correspondence 
between the set of ideals of $\cpi$ and the closed subset of $\Prim(\cpi)$.
By Proposition \ref{closed}, 
the closed subset of $\Prim(\cpi)$ corresponds bijectively 
to the set of $\omega$-invariant pairs.
\eprf

\section{More about $\cpi$}

In this section, we gather some general results on $\cpi$.
First we compute the strong Connes spectrum of the action 
$\alpha^{\omega}:G\curvearrowright\Oi$.
We need the following lemma.

\blem\label{H_csg}
For any $\omega\in\Gamma^\infty$, we have 
$\{0\}\cup H_{\csg}
=\{0\}\cup\bigcap_{n=1}^\infty\overline{\bigcup_{i=n}^\infty(\sg+\omega_i)}$.
\elem

\bprf
It suffices to show that 
$$\csg\setminus\big(\{0\}\cup\bigcup_{i=1}^\infty(\csg+\omega_i)\big)\subset 
\overline{\bigcup_{i=n+1}^\infty(\sg+\omega_i)}$$
for any $n\in\Z_+$.
Take $\gamma\in\csg\setminus
\big(\{0\}\cup\bigcup_{i=1}^\infty(\csg+\omega_i)\big)$ 
and $n\in\Z_+$.
Since $\gamma\in\csg$, there exists a sequence $\{\mu_k\}\subset\W_\infty$
such that $\gamma=\lim_{k\to\infty}\omega_{\mu_k}$.
We will show that we can find an integer grater than $n$ in the word $\mu_k$
for infinitely many $k$, from which it follows that 
$\gamma\in\overline{\bigcup_{i=n+1}^\infty(\sg+\omega_i)}$.
To the contrary, assume that $\mu_k\in\W_n$ for sufficiently large $k$.
Then there exists $i\in\{1,2,\ldots,n\}$ which appears in $\mu_k$ eventually.
We have $\gamma-\omega_i=\lim_{k\to\infty}(\omega_{\mu_k}-\omega_i)\in\csg$.
This contradicts the fact that $\gamma\notin\csg+\omega_i$.
Hence 
$\{0\}\cup H_{\csg}
=\{0\}\cup\bigcap_{n=1}^\infty\overline{\bigcup_{i=n}^\infty(\sg+\omega_i)}$.
\eprf

\bpr\label{SCS}
The strong Connes spectrum $\widetilde{\Gamma}(\alpha^{\omega})$ 
of the action $\alpha^{\omega}$ is $\{0\}\cup H_{\csg}$.
\epr

\bprf
By \cite[Lemma 3.4]{Ki}, we have 
$$\widetilde{\Gamma}(\alpha^{\omega})
=\{\gamma\in\Gamma\mid \widehat{\alpha^{\omega}}_\gamma(I)\subset I, 
\mbox{for any ideal } I \mbox{ of } \cpi\},$$
where $\widehat{\alpha^{\omega}}:\Gamma\curvearrowright\cpi$ is the dual
action of $\alpha^{\omega}$.
For an $\omega$-invariant pair $\widetilde{X}=(X,X^{(\infty)})$ 
and $\gamma\in\Gamma$, 
we see that $\widehat{\alpha^{\omega}}_\gamma(I_{\widetilde{X}})
=I_{\widetilde{X}-\gamma}$ where 
$\widetilde{X}-\gamma=(X-\gamma,X^{(\infty)}-\gamma)$.
Hence $\widehat{\alpha^{\omega}}_\gamma(I_{\widetilde{X}})\subset 
I_{\widetilde{X}}$ is equivalent 
to say that $X+\gamma\subset X$ and $X^{(\infty)}+\gamma\subset X^{(\infty)}$ 
for an $\omega$-invariant pair $\widetilde{X}=(X,X^{(\infty)})$ 
and $\gamma\in\Gamma$.
Considering the case that $\widetilde{X}=(\csg,\{0\}\cup H_{\csg})$,
we have $(\{0\}\cup H_{\csg})+\gamma\subset\{0\}\cup H_{\csg}$
for $\gamma\in\widetilde{\Gamma}(\alpha^{\omega})$.
Hence $\widetilde{\Gamma}(\alpha^{\omega})\subset\{0\}\cup H_{\csg}$.
Let $(X,X^{(\infty)})$ be an $\omega$-invariant pair. 
For $\gamma\in X$, we get 
$$\gamma+\bigcap_{n=1}^\infty\overline{\bigcup_{i=n}^\infty(\sg+\omega_i)}
=\bigcap_{n=1}^\infty\overline{\bigcup_{i=n}^\infty(\gamma+\sg+\omega_i)}
\subset\bigcap_{n=1}^\infty\overline{\bigcup_{i=n}^\infty(X+\omega_i)}
\subset H_X\subset X^{(\infty)}.$$
By Lemma \ref{H_csg}, we have 
$X^{(\infty)}+(\{0\}\cup H_{\csg})\subset X^{(\infty)}$.
Since $\{0\}\cup H_{\csg}\subset\csg$, 
we have $X+(\{0\}\cup H_{\csg})\subset X$.
Hence when $\omega$ satisfies Condition \ref{cond}, we have
$\widetilde{\Gamma}(\alpha^{\omega})\supset \{0\}\cup H_{\csg}$
by Theorem \ref{idestr1}, 
and so $\widetilde{\Gamma}(\alpha^{\omega})=\{0\}\cup H_{\csg}$.
Next we consider the case 
that $\omega$ does not satisfy Condition \ref{cond}.
For an $\omega$-invariant pair $(Y,X^{(\infty)})$, 
we have $X+(H_{\csg}\setminus\{0\})\subset X^{(\infty)}$ 
by the former part of this proof, where 
$X=\{\gamma\in\Gamma\mid 
([\gamma],\theta)\in Y \mbox{ for some }\theta\in\T\}$.
Hence for any $([\gamma_0],\theta_0)\in Y$ 
and $\gamma\in H_{\csg}\setminus\{0\}$,
we have $\gamma_0+\gamma\in X^{(\infty)}$ 
because $\gamma_0\in X$.
Since $[X^{(\infty)}]\times\T\subset Y$, we have 
$([\gamma_0+\gamma],\theta_0)\in Y$. 
Therefore we also have 
$\{0\}\cup H_{\csg}\subset\widetilde{\Gamma}(\alpha^{\omega})$
by Theorem \ref{idestr2}.
Thus $\widetilde{\Gamma}(\alpha^{\omega})=\{0\}\cup H_{\csg}$.
\eprf

Next we give necessary and sufficient conditions for $\omega\in\Gamma^\infty$ 
that the crossed product $\cpi$ becomes simple or primitive.

\blem\label{0}
Let $I$ be an ideal of the crossed product $\cpi$. 
Then $I=0$ if and only if $X_I=\Gamma$. 
\elem

\bprf
The ``only if'' part is trivial. 
One can easily prove the ``if'' part 
by the same arguments as in the proofs of 
Proposition \ref{cond.exp2} and Theorem \ref{idestr1}.
\eprf

\bpr
For $\omega\in\Gamma^\infty$, the following are equivalent:
\benu
\item The crossed product $\cpi$ is simple.
\item There are no $\omega$-invariants sets 
other than $\Gamma$ and $\emptyset$.
\item $\Gamma=\csg$.
\eenu

If $\cpi$ is simple, then it is purely infinite.
\epr

\bprf
The equivalence between (i) and (ii) follows from Lemma \ref{0}.
(ii) implies (iii) because $\csg$ is $\omega$-invariant. 
(iii) implies (ii) because $X=X+\csg$ if $X$ is $\omega$-invariant.

For the last statement, see \cite[Proposition 5.2]{Ka2}.
\eprf

The equivalence between (i) and (iii) was already proved 
by A. Kishimoto \cite{Ki} by using strong Connes spectrum.
Note that the strong Connes spectrum $\widetilde{\Gamma}(\alpha^{\omega})$
is equal to $\Gamma$ if and only if $\csg=\Gamma$ by Proposition \ref{SCS}.

\bpr
The following conditions for $\omega\in\Gamma^\infty$ are equivalent:
\benu
\item The crossed product $\cpi$ is primitive.
\item $\Gamma$ is a prime $\omega$-invariant set.
\item The closed group generated by $\omega_1,\omega_2,\ldots$ 
is equal to $\Gamma$.
\eenu
\epr

\bprf
(i)$\Rightarrow$(ii): This follows from Proposition \ref{prime}.

(ii)$\Rightarrow$(i): It suffices to show that $0$ is prime.
Let $I_1,I_2$ be ideals of $\cpi$ with $I_1\cap I_2=0$.
We have $X_{I_1}\cup X_{I_2}=X_{I_1\cap I_2}=\Gamma$.
Since $\Gamma$ is prime, either $X_{I_1}\supset\Gamma$ or $X_{I_2}\supset\Gamma$.
If $X_{I_1}\supset\Gamma$ hence $X_{I_1}=\Gamma$, then $I_1=0$ by Lemma \ref{0}.
Similarly if $X_{I_2}\supset\Gamma$, then $I_2=0$.
Thus $0$ is prime and so $\cpi$ is a primitive $C^*$-algebra.

(ii)$\iff$(iii): This follows from Proposition \ref{Xprime}.
\eprf

One can prove the equivalence between (i) and (iii) in the above theorem
by characterization of primitivity of crossed products in terms of 
the Connes spectrum due to D. Olesen and G. K. Pedersen \cite{OP} and
the computation of the Connes spectrum of our actions $\alpha^{\omega}$ 
due to A. Kishimoto \cite{Ki}.

\bpr\label{CP}
The crossed product $\cpi$ is isomorphic to 
the Cuntz-Pimsner algebra ${\mathcal O}_E$ of $C_0(\Gamma)$-bimodule
$E=C_0(\Gamma)^\infty$, whose left module structure is given by 
$$f\cdot(f_1,f_2,\ldots,f_n,\ldots)=(\sigma_{\omega_1}(f)f_1,\sigma_{\omega_2}(f)f_2,\ldots,\sigma_{\omega_n}(f)f_n,\ldots)\in E$$
for $f\in C_0(\Gamma)$ and $(f_1,f_2,\ldots,f_n,\ldots)\in E$.
\epr

\bprf
The inclusion $C_0(\Gamma)\hookrightarrow\cpi$ and 
$E\ni (0,\ldots,0,f_n,0\ldots)\mapsto S_nf_n\in\cpi$ 
satisfies the conditions in \cite[Theorem 3.12]{Pi}. 
Hence there exists a $*$-homomorphism 
$\varphi:{\mathcal O}_E\to\cpi$ which is surjective 
since $\cpi$ is generated by $\{S_nf\mid n\in\Z_+,\ f\in C_0(\Gamma)\}$.
One can show that $\varphi$ is injective by using Lemma \ref{isom}.
Thus $\cpi$ is isomorphic to ${\mathcal O}_E$.
\eprf

\bco
The inclusion $C_0(\Gamma)\hookrightarrow\cpi$ is a KK-equivalence.
Hence for $i=0,1$, we have $K_i(\cpi)=K_i(C_0(\Gamma))$.
\eco

\bprf
See \cite[Corollary 4.5]{Pi}.
\eprf

\bpr\label{embedinfty}
If $\omega\in\Gamma^\infty$ satisfies 
$-\omega_i\notin\overline{\{\omega_{\mu}\mid\mu\in\W_n\}}$ 
for any $i,n\in\Z_+$, 
then the crossed product $\cpi$ is AF-embeddable.
\epr

\bprf
See \cite[Proposition 5.1]{Ka2}.
\eprf

\end{document}